\documentclass[12pt,reqno]{amsart}
\usepackage{amsmath, amsfonts, amssymb, amsthm, amscd, amsbsy}
\usepackage{fancyhdr}
\usepackage[usenames,dvipsnames,svgnames,x11names,hyperref]{xcolor}
\usepackage{geometry}
\usepackage{graphicx}
\usepackage{hyperref}
\usepackage{verbatim}
\usepackage{dsfont}
\usepackage{enumitem}

\numberwithin{equation}{section}
\providecommand{\U}[1]{\protect\rule{.1in}{.1in}}

\hypersetup{
colorlinks=true,
breaklinks=true,
urlcolor=NavyBlue,
linkcolor=Fuchsia,
bookmarksopen=false,
filecolor=black,
citecolor=ForestGreen,
linkbordercolor=red
}
\allowdisplaybreaks

\newtheorem{theorem}{Theorem}[section]

\newtheorem{prop}[theorem]{Proposition}

\newtheorem{lemma}[theorem]{Lemma}
\newtheorem{remark}[theorem]{Remark}
\newtheorem{example}[theorem]{Example}
\newtheorem{examples}[theorem]{Examples}
\newtheorem{foo}[theorem]{Remarks}

\def\vint{\mathop{\mathchoice%
          {\setbox0\hbox{$\displaystyle\intop$}\kern 0.22\wd0%
           \vcenter{\hrule width 0.6\wd0}\kern -0.82\wd0}%
          {\setbox0\hbox{$\textstyle\intop$}\kern 0.2\wd0%
           \vcenter{\hrule width 0.6\wd0}\kern -0.8\wd0}%
          {\setbox0\hbox{$\scriptstyle\intop$}\kern 0.2\wd0%
           \vcenter{\hrule width 0.6\wd0}\kern -0.8\wd0}%
          {\setbox0\hbox{$\scriptscriptstyle\intop$}\kern 0.2\wd0%
           \vcenter{\hrule width 0.6\wd0}\kern -0.8\wd0}}%
          \mathopen{}\int}

\newcommand{\hs}{\mathbb{H}^n}
\newcommand{\bn}{\mathbb{B}^n}
\newcommand{\rn}{\mathbb{R}^n}
\newcommand{\halfrn}{\mathbb{R}^n_+}
\newcommand{\inth}{\int_{\mathbb{H}^n}}
\newcommand{\cs}{C^\infty_c}

\begin{document}

\title{Symmetry of solutions to higher and fractional order semilinear equations on hyperbolic spaces}
\thanks{The second and third authors' research is partly supported by a Simons Collaboration grant from the Simons Foundations.}
\author{Jungang Li, Guozhen Lu, Jianxiong Wang}
\address{Jungang Li: Department of Mathematics\\Brown University\\Providence, RI 02912, USA.}
\email{jungang\_li@brown.edu}
\address{Guozhen Lu: Department of Mathematics\\University of Connecticut\\Storrs, CT 06269, USA.}
\email{guozhen.lu@uconn.edu}
\address{Jianxiong Wang: Department of Mathematics\\University of Connecticut\\Storrs, CT 06269, USA.}
\email{jianxiong.wang@uconn.edu}
\date{}
\setlength\parindent{0pt}

\begin{abstract}
    We show that nontrivial solutions to higher and fractional order equations with certain nonlinearity are radially symmetric and nonincreasing
    on geodesic balls in the hyperbolic space $\hs$ as well as on the entire space $\hs$. Applying the Helgason-Fourier analysis techniques on $\hs$, we
    develop a moving plane approach for integral equations on $\hs$.
     We also establish the symmetry to solutions of certain equations with singular terms on Euclidean spaces.
    Moreover, we obtain symmetry to solutions of some semilinear equations involving fractional order derivatives.
\end{abstract}
\maketitle
\smallskip

\section{Introduction}

In the celebrated work of Gidas, Ni and Nirenberg \cite{GNN}, they considered the following boundary value problem on the ball $B_R(0)\subset{\rn}$.
\begin{equation*}
    \begin{cases}
        -\Delta u = f(u) &\text{ in } B_R(0)\\
         u=0 &\text{ on } \partial B_R(0),
    \end{cases}
\end{equation*}
where $f$ is of class $C^1$. They proved that any positive solution $u$ in $C^2(\overline{B_R(0)})$
is radially symmetric and decreasing. Their approach is the so-called moving plane method, which was initiated by the Soviet mathematician Alexandrov in the 1950s, and further developed by Serrin \cite{serrin}.
Subsequently, Caffarelli, Gidas and Spruck \cite{Caff} proved the symmetry of solutions when $f$ is with critical growth. To be more precise, they considered the nonlinearity term $f$ which is a locally nondecreasing Lipschitz function such that $f(0)=0$, and satisfies the following growth condition: for sufficiently large $t$, the function
 $t^{-\frac{n+2}{n-1}}f(t)$ is nonincreasing and
 $f(t)\geq ct^p$ for some $p\geq \frac{n}{n-2}$. In the past decades, the moving plane method together with its applications to nonlinear PDEs have been extensively studied.
We refer the interested reader to, e.g., \cite{Bere,  W-X, ChenLi1} and the vast list of references therein. In particular, the moving plane method in integral form in Euclidean spaces has been developed in \cite{CLO} which motivated our approach in the hyperbolic space $\mathbb{H}^n$ in the current paper.
\medskip

The moving plane method has also been developed on non-Euclidean spaces. In \cite{K-P2}, Kumaresan and Prajapat established a Serrin type symmetry result on bounded domains in $\hs$.
In a subsequent article \cite{K-P1}, the same authors obtained the analogous result of Gidas-Ni-Nirenberg type on hyperbolic spaces and spheres. For instance, they proved:

\medskip

{\bf Theorem A.}
{\it
  Let $\Omega$ be a geodesic ball in $\hs$ and $u\in C^2(\bar{\Omega})$ be a positive solution of the equation
  \begin{equation*}
    \begin{cases}
      -\Delta_{\mathbb{H}^n} u = f(u) &\text{ on } \Omega, \\
       u=0  &\text{ on } \partial\Omega,
    \end{cases}
  \end{equation*}
  where $\Delta_{\mathbb{H}^n}$ is the Laplace-Beltrami operator on $\hs$, $f$ is a $C^1$ function. Then, $u$ is radially symmetric.}

  \medskip

As for the entire space $\hs$, Almeida, Damascelli and Ge \cite{ADG} studied the following problem:
\begin{equation}\label{eq1.1}
    \begin{cases}
    -\Delta_{\mathbb{H}^n}u=f(u) &\text{ on } \hs, \\
    u>0 &\text{ on } \hs,
\end{cases}
\end{equation}
where $f: (0, \infty) \to \mathbb{R}$ is locally Lipschitz continuous.
Among other results, they proved the following:

\medskip
{\bf Theorem B.} {\it
Let $u\in C^1(\hs)$ be a weak solution to the problem \eqref{eq1.1}, then $u$ is radially symmetric and decreasing with respect to some point $P \in \mathbb{H}^n$, provided one of the following conditions:
\begin{enumerate}
    \item $u \to 0$ as $|x| \to \infty$ and there exists $s_0 > 0$ such that $f$ is nonincreasing on $(0 , s_0)$.

    \item $u \to 0$ as $|x| \to \infty$ and there exists $s_0,\alpha > 0$ such that if $0 < a,b < s_0$, $\frac{f(a) - f(b)}{a - b} \leq G + (a + b)^\alpha$ and $u \in L^2 (\mathbb{H}^n) \cap H^{\alpha n /2}(\hs)$, where $G, C$ are some dimensional positive constants.

    \item $u \in W^{1,2}(\mathbb{H}^n)$, there exists $\alpha > 0 $ such that for $0 < a< b$, $\left| \frac{f(a) - f(b)}{a - b} \right| \leq G + C(a + b )^\alpha$ and $u \in L^2 (\mathbb{H}^n) \cap H^{\alpha n /2}(\hs)$.
\end{enumerate}
}
\medskip

The main purpose of the present paper is to investigate the symmetry of solutions of higher order and fractional order semilinear equations with critical growth on hyperbolic spaces. To illustrate our main results, we would like to first recall some previous results on semilinear equations with critical growth. In the celebrated work \cite{B-N}, the following so-called Br\'ezis-Nirenberg problem was systematically studied:
\begin{equation}\label{B-N}
    \begin{cases}
        -\Delta u - \lambda u= u^{2^*-1}  &\text{on } \Omega\\
        u>0 &\text{on } \Omega\\
        u=0 &\text{on } \partial\Omega,
    \end{cases}
\end{equation}
where $\Omega$ is a bounded domain in $\mathbb{R}^n \text{ for } n\geq 3$ and $2^*=\frac{2n}{n-2}$, which is the critical Sobolev exponent. Br\'ezis and Nirenberg proved the following result:

\medskip
{\bf Theorem C.}
{\it
When $n\geq 4$, probelm (\ref{B-N})
has a nontrivial solution for every $\lambda\in(0,\Lambda_1 (\Omega))$, where $\Lambda_1 (\Omega)$ is the first Dirichlet eigenvalue of $-\Delta$. Morever, when $\Omega$ is a star-shaped domanin, the problem (\ref{B-N}) has no solution if $\lambda\leq 0$.
In particular, when $\Omega=\mathbb{B}^3 \subset \mathbb{R}^3$,  (\ref{B-N}) has nontrivial solution
if and only if $\lambda\in(\frac{1}{4}\Lambda_1 (\Omega),\Lambda_1 (\Omega))$.
}

\medskip

Due to the lack of compactness caused by the appearance of the critical Sobolev exponent, the Br\'ezis-Nirenberg problem cannot be studied directly by the variational method. Such a phenomenon occurs in many differential geometry and PDE problems. One typical example is the famous Yamabe problem (we refer the interested readers to, e.g., \cite{Yamabe1,LeeParker1,Yamabe2,Yamabe3,Yamabe} and the references therein). After \cite{B-N}, different versions of higher order Br\'ezis-Nirenberg problems have been extensively studied in the past decades. One of them can be formulated as follows:
\begin{equation}\label{higher B-N}
    \begin{cases}
        (-\Delta)^k u = \lambda u+ |u|^{q-1}u &\text{on } \Omega\\
        u=\nabla u=\cdots=\nabla^{k-1}u=0 &\text{on } \partial\Omega,
    \end{cases}
\end{equation}
where $\Omega\subset\mathbb{R}^n$ is a bounded domain, $n>2k$ and $q=\frac{n+2k}{n-2k}$ is the corresponding critical Sobolev exponent. Several existence results analogous to the second order
Br\'ezis-Nirenberg problem have been obtained. For example, Gazzola \cite{Gazzola} established the following:

\medskip
{\bf Theorem D.}
{\it
Let $\Lambda_1((-\Delta)^k,\Omega)$ denote the first Dirichlet eigenvalue of $(-\Delta)^k$ on $\Omega$, then
\begin{enumerate}
    \item When $n\geq 4k$, there exists a solution $u\in W^{k,2}_0(\Omega)$ to the Dirichlet problem (\ref{higher B-N}) for every $\lambda\in (0,\Lambda_1((-\Delta)^k,\Omega))$;
    \item When $2k+1\leq n\leq 4k-1$, there exists $0<\bar{\Lambda}<\Lambda_1((-\Delta)^k,\Omega)$ such that for every $\lambda\in (\bar{\Lambda},\Lambda_1((-\Delta)^k,\Omega))$, (\ref{higher B-N}) has a solution $u\in W^{k,2}_0(\Omega)$.
\end{enumerate}
}
\medskip

In particular, when $\Omega=B^n$, a Euclidean ball, Grunau \cite{Grunau} obtained stronger results. It has been shown that the solution of (\ref{higher B-N}) actually belongs to $C^\infty(B^n)\cap C^{2k+1}(\overline{B^n})$
and the solution is positive, radially symmetric and decreasing.  When $\Omega$ is a ball, Pucci and Serrin \cite{PucciSerrin1} conjectured that for dimensions $n = 2k + 1, 2k + 2, \cdots, 4k - 1$, the necessary condition for the existence of solution is that $\lambda$ should be larger than some positive constant number. Such dimensions are called critical dimensions. When $n = 3$ and $k = 1$, Brezis and Nirenberg \cite{B-N} already found such a lower bound explicitly. Pucci and Serrin proved that $n = 2k + 1$ is critical. When $k = 2$, the biharmonic version of Pucci-Serrin's conjecture has been proved by Edmunds, Fortunato and Jannelli \cite{EdmundsFortunatoJannelli1}. The cases $k = 3,4$ were due to Bernis and Grunau \cite{BernisGrunau1, Grunau}.
\\~\\
Due to the importance in the study of differential geometry, Br\'ezis-Nirenberg problem also receives much attention in non-Euclidean settings.
One typical example is the hypberbolic space $\hs$, which is a complete noncompact simply connected Riemannian manifold with constant sectional curvature equals $-1$.
%Among many different models for hypberbolic spaces, we first recall the Poincar\'e ball model, which is given by the open unit ball $\mathbb{B}^n\subset\rn$ equipped with the Poincar\'e metric %
%$ds^2=\frac{4\left(dx_1^2+\cdots+dx_n^2\right)}{\left(1-|x|^2\right)^2}%$.
Mancini and Sandeep \cite{Mancini-Sandeep} studied the following second order semilinear equation:
$$-\Delta_{\hs}u=\lambda u+u^q,$$
where $\lambda$ is a real parameter, $q>1$ if $n=2$ and $1<q\leq \frac{n+2}{n-2}$ if $n \geq 3$.
They obtained several existence/nonexistence results, depending on the range of $\lambda$ in different dimensions.
Moreover, they proved that the positve solution has hypberbolic symmetry, if $\lambda\leq \left(\frac{n-1}{2}\right)^2, n\geq 2$, i.e., there exists $x_0\in\hs$ such that
$u$ is constant on hypberbolic spheres centered at $x_0$.
\\~\\
Recently, Li, Lu and Yang \cite{LLY} studied the higher order Br\'ezis-Nirenberg problem associated with GJMS operators on hypberbolic spaces.
The GJMS operator $P_k$ is of a class of 2k-th order conformal operators and has been wildly studied in the conformal geometry, see e.g. \cite{FeffermanGr,FeffermanGr1,GJMS2,GJMS3}.
More specifically, the GJMS operator on $\mathbb{H}^n$ can be inductively defined as follows:
\begin{equation*}
    P_k=P_1(P_1+2)\cdots(P_1+k(k-1)), k\in\mathbb{N},
\end{equation*}
where $P_1=-\Delta_{\mathbb{H}^n}-\frac{n(n-2)}{4}$ is the conformal Laplacian on $\mathbb{H}^n$ and
\begin{equation*}
    P_2=\left(-\Delta_{\mathbb{H}^n}-\frac{n(n-2)}{4}\right)\left(-\Delta_{\mathbb{H}^n}-\frac{(n+2)(n-4)}{4}\right)
\end{equation*}
is the Paneitz operator. We use $W^{k,2}_0(\Omega)$ to denote the Sobolev space on the hyperbolic space $\mathbb{H}^n$ as defined in Subsection \ref{Subsection2.2}. In \cite{LLY}, among other results, the authors proved the following:

\medskip
{\bf Theorem E.}
{\it
Consider
\begin{equation}\label{PDE0}
    \begin{cases}
        P_k u -\lambda u = |u|^{q-2}u \;&\text{ in } \Omega,\\
        \nabla_{\hs}^\alpha u = 0,\; \alpha=0,\cdots,k-1 &\text{ on } \partial\Omega,
    \end{cases}
\end{equation}
where $q=\frac{2n}{n-2k}, 2\leq k<\frac{n}{2}, \Omega$ is a bounded domain in $\hs$ with $C^1$ boundary and $\nabla_{\mathbb{H}^n}$ is the hypberbolic gradient (which will be defiend in Section 2).
Let $\Lambda_1(P_k, \Omega)$ be the first Dirichlet eigenvalue of $P_k$ on $\Omega$ defined by
$$\Lambda_1(P_k,\Omega)=\inf_{u\in C^\infty_0\setminus\{0\}}\frac{\int_\Omega P_ku\cdot udV}{\int_\Omega |u|^2dV},$$
where $dV$ is the volume element in $\hs$. Then the following holds:
\begin{enumerate}
    \item when $n\geq 4k $ and $0<\Lambda<\Lambda_1 (P_k,\Omega)$, (\ref{PDE0}) has at least one nontrivial solution in $W^{k,2}_0(\Omega)$;
    \item when $2k+1\leq n\leq 4k-1$, there exists a positive constant $\Lambda^*$ such that for $\Lambda^*<\lambda<\Lambda_1 (P_k,\Omega)$, (\ref{PDE0}) has at least one nontrivial solution in $W^{k,2}_0(\Omega)$.
\end{enumerate}
}
\medskip

  It is worth pointing out that in \cite{LLY}, the authors also investigated the case when $\Omega = \mathbb{H}^n$ and no boundary condition is assigned. Their result can be viewed as a higher order generalization of \cite{Mancini-Sandeep}. When $k = 1$, (\ref{PDE0}) was initiated in \cite{Stapelkamp} and further studied in \cite{Benguria}. Unlike the second order case, due to the complexity of the operator, one cannot expect to reduce higher order Br\'ezis-Nirenberg problems on $\mathbb{H}^n$ to the corresponding ones on Euclidean spaces and apply known results from there. The authors in \cite{LLY} developed a new approach, which combines the knowledge of higher order Hardy-Sobolev-Maz'ya's inequalities, Helgason-Fourier analysis and Green's function estimate on hyperbolic spaces. These ingredients will still play a central role in the present paper. \\

 On the other hand, the symmetry of solutions plays an important role in the study of Br\'ezis-Nirenberg problems, especially in lower dimension cases. This is not surprising since the symmetry (or the invariance under group actions) is the key ingredient to the characterization of extremal functions of sharp Sobolev type inequalities (see e.g. \cite{B-N,GNN,Talenti1}) and it has been shown in \cite{B-N} that the solvability of Br\'ezis-Nirenberg problem is deeply related to the sharp Sobolev inequality. Nevertheless, the classical moving plane method heavily depends on the maximum principle, which is not applicable to higher order equations.  In the work of Chen, Li and Ou \cite{CLO}, a new moving plane method in terms of integral equations was developed and applied to study higher order differential equations and integral equations (see also \cite[Chapter 8]{ChenLi1}). Subsequently, Lu and Zhu \cite{L-Z} studied the moving plane method on upper half spaces. Inspired by these works, with the help of Hardy-Littlewood-Sobolev inequalities and the Helgason-Fourier analysis on $\mathbb{H}^n$, the authors in \cite{LLY} developed a moving plane method to integral forms on $\mathbb{H}^n$ and applied it to study the symmetry of solutions to the higher order Br\'ezis-Nirenberg problems on $\mathbb{H}^n$. To be precise, they established the following:

\medskip
{\bf Theorem F.}
{\it
Let $k\geq 2$. If $u\in W^{k,2}_0(\hs)$ is a positive weak solution of the equation
$$P_ku-\lambda u = |u|^{q-2}u,$$
then there exists a point $P \in\hs$ such that $u$ is constant on the geodesic spheres centered at $P$. Morever, $u$ is non-increasing.
}

\medskip

  One of our primary goals of the present paper is to investigate the symmetry of solutions when the  equation is equipped with the generic forcing term:
\begin{equation}\label{P_k}
        P_k u = f(u)\text{ on } \hs.
\end{equation}
Our first result is the following:

\begin{theorem}\label{thm1}
    Let $k\geq 2, f$ be Lipschitz continuous, non-decreasing. Assume $f^\prime(u)\in L^{\frac{n}{2k}}(\hs)$. If $u\in W^{k,2}_0(\mathbb{H}^n)$ is a positive solution of the equation (\ref{P_k}),
    then there exists a point $ P \in \mathbb{H}^n$ such that u is constant on the geodesic spheres centered
    at $P$. Moreover, $u$ is non-increasing.
\end{theorem}

  The moving plane method to integral forms requires certain monotonicity of the corresponding Green's function. In \cite{LLY}, the monotonicity was obtained by a fraction of the operator $P_k - \lambda$ and a complicated analysis to the Green's function. On the other hand, an explicit expression of Green's function $(-\Delta_{\hs}+\lambda)^{-1}$ for $\lambda>-\frac{(n-1)^2}{4}$ is given by
(see for \cite{HM} $\lambda\geq 0$ and \cite{Lhq} for $\lambda>-\frac{(n-1)^2}{4}$)
\begin{equation}\label{Gf_1}
    (-\Delta_{\hs}+\lambda)^{-1}=(2\pi)^{-\frac{n}{2}}(\sinh\rho)^{-\frac{n-2}{2}}e^{-\frac{(n-2)\pi}{2}i}Q^{\frac{n-2}{2}}_{\theta_n(\lambda)}(\cosh\rho),
\end{equation}
where $\theta_n(\lambda)=\sqrt{\lambda+\frac{(n-1)^2}{4}}-\frac{1}{2}$ and $Q^\mu_\nu(z)$ is the Legendre function of second type which will be defined in the next lemma.
 Lu and Yang \cite{LY} recently proved that (\ref{Gf_1}) is also vaild for $\lambda=-\frac{(n-1)^2}{4}$. More precisely,
\begin{lemma}
    Let $n\geq 3$, we have
    \begin{align*}
        \left(-\Delta_{\mathbb{H}^n}-\frac{(n-1)^2}{4}\right)^{-1}=(2\pi)^{-\frac{n}{2}}(\sinh\rho)^{-\frac{n-2}{2}}e^{-\frac{(n-2)\pi}{2}i}Q^{\frac{n-2}{2}}_{-\frac{1}{2}}(\cosh\rho),
    \end{align*}
    where $Q_\nu^\mu(z)$ is the Legendre function of second type defiend by
    \begin{align*}
        Q_\nu^\mu(z)=&e^{i(\pi\mu)}2^{-\nu-1}\frac{\Gamma(\nu+\mu+1)}{\Gamma(\nu+1)}(z^2-1)^{-\mu/2}\int_0^{\pi}(z+\cos t)^{\mu-\nu-1}(\sin t)^{2\nu+1}dt,\\
        &\text{Re}\nu>1, \text{Re}(\nu+\mu+1)>0.
    \end{align*}
\end{lemma}

   We realize that the hypergeometric function expression to the Green's function will directly give us the monotonicity and hence significantly simplify our argument to the equation (\ref{P_k}).

\medskip

Our second main result gives the symmetry of solutions to the following Dirichlet problem:
\begin{equation}\label{P_k*}
    \begin{cases}
        P_k u = f(u) &\text{ on } B,\\
        \nabla_{\hs}^\alpha u = 0,\; \alpha \leq k-1 &\text{ on }\partial B,
    \end{cases}
\end{equation}
where $B = B_{P}(\rho)$ is a geodesic ball centered at $P \in \mathbb{H}^n$ with radius $\rho$. Our result reads as follows:

\begin{theorem}\label{thm2}
    Assume that $f: [0,\infty)\rightarrow\mathbb{R}$ is a continuous, non-decreasing function with $f(0)\geq 0$
    and that $u\in W^{k,2}_0(B)\cap L^\infty(B)$ is a positive weak solution to the Dirichlet problem (\ref{P_k*}).
    Then $u$ is radially symmetric and strictly decreasing in the radial variable.
\end{theorem}

We would like to add some comments on Theorem \ref{thm2}. As we mentioned previously, the classical moving plane method fails on higher order equations, due to the lack of maximum principle.
To overcome this difficulty, Berchio, Gazzola and Weth \cite{BGW} proved an alternative Hopf lemma towards the polyharmonic Dirichlet problem on $\mathbb{R}^n$ and hence re-developed the moving plane method. For instance, they proved the following:

\medskip
{\bf Theorem G.}
{\it
  Assume that $f : [0, \infty) \to \mathbb{R}$ is a continuous nondecreasing function with $f(0) \geq 0$ and $u \in W^{k,2}_0 \cap L^\infty$ defined on the Euclidean ball with radius $R$, i.e. $B = B_R$, is a nonnegative nontrivial weak solution to the following Dirichlet problem

  $$
    \begin{cases}
      (- \Delta)^k u = f(u) &\text{ on } B, \\
      \nabla^\alpha u = 0, \alpha = 1, \cdots , k-1 &\text{ on } \partial B.
    \end{cases}
  $$
  Then $u$ is radially symmetric and strictly decreasing with respect to the radial variable.
}

\medskip

Theorem \ref{thm2} can be viewed as a hyperbolic version of the above Theorem G. Our main contribution in Theorem \ref{thm2} is that we obtain a Boggio type formula for the corresponding Green's function and further use it to establish a Hopf type lemma on hyperbolic spaces and this enables us to perform the moving plane argument on $\mathbb{H}^n$. As far as we know, Theorem \ref{thm2} is the first moving plane result for higher order equations with general nonlinearity $f(u)$  on bounded domains in $\mathbb{H}^n$. Moreover, Berchio-Gazzola-Weth's theorem is valid if $f(u)$ is replaced by the nonautonomous radial nonlinearity $f(|x|, u)$ provided that $f$ is continuous and nonincreasing with respect to the first variable. We realize that Theorem \ref{thm2} will imply Berchio-Gazzola-Weth's theorem when the monotonicity of $f$ breaks. To be precise, we first recall the Poincar\'e ball model of $\mathbb{H}^n$, which is the unit ball in $\mathbb{R}^n$ equipped with the Poincar\'e metric $ds^2 = \frac{4 (dx_1^2 + \cdots | dx_n^2)}{(1 - |x|^2)^2}$ (see Section 2 for more details). Due to the invariance of the operator $P_k$ under the conformal transform, without loss of generality, one only needs to prove Theorem \ref{thm2} when $P$ is the origin. It is well known in the hyperbolic geometry that the geodesic ball centered at the origin coincides with the Euclidean ball (with different radius though). On the other hand, Liu \cite{Liu} proved the following identity for any $u\in C^\infty_0(\hs)$ (in the Poicar\'e ball model):
$$P_ku=\left(\frac{1-|x|^2}{2}\right)^{k+\frac{n}{2}}(-\Delta)^k\left((\frac{1-|x|^2}{2})^{k-\frac{n}{2}}u\right).$$
Therefore, (\ref{P_k}) can be written as an equivalent equation on the Euclidean ball $B_1(0)$:
\begin{equation}\label{eq1.7}
    (-\Delta)^kv=\left(\frac{2}{1-|x|^2}\right)^{k+\frac{n}{2}}f((\frac{1-|x|^2}{2})^{\frac{n}{2}-k}v).
\end{equation}
Thus, Theorem \ref{thm2} implies the following:
\begin{theorem}\label{thm3}
    Let $B=B(0,R^\prime)$ be the Euclidean ball centered at the origin with radius $R^\prime$.
    Assume $f$ satisfies the same conditions as in Theorem \ref{thm2}.
    If $v\in W^{k,2}_0(B)$ is a positive weak solution of the equation:
    \begin{equation}
        \begin{cases}
            (-\Delta)^kv=\left(\frac{2}{1-|x|^2}\right)^{k+\frac{n}{2}}f((\frac{1-|x|^2}{2})^{\frac{n}{2}-k}v)&\text{ in }B(0,R^\prime),\\
            \nabla^\alpha v=0, \quad \alpha = 1, \cdots , k-1 &\text{ on } \partial B(0,R^\prime),
        \end{cases}
    \end{equation}
    then $v$ is radially symmetric with respect to the origin and strctly decreasing.
\end{theorem}

\begin{remark}
Notice that the nonlinearity term
\begin{equation}\label{nonauto_g}
g(|x|,v)=\left(\frac{2}{1-|x|^2}\right)^{k+\frac{n}{2}}f((\frac{1-|x|^2}{2})^{\frac{n}{2}-k}v)
\end{equation}
does not satisfy the monotonicity conditions in Berchio-Gazzola-Weth's thoerem \cite{BGW}, and Theorem \ref{thm2} actually broadens the class of nonlinearity in \cite{BGW}.
\end{remark}

We now recall the upper half space model of the hypberbolic space, which is the upper half space $\mathbb{R}^n_+ = \{ (x_1, \cdots , x_n) \in \mathbb{R}^n :  x_n > 0\}$ equipped with the metric $ds^2=\frac{dx_1^2+\cdots+dx_n^2}{x_n^2}$.

In \cite{L-Z}, Lu and Zhu studied the following integral equation on upper half spaces.
\begin{equation}
    u(x)=\int_{\halfrn}G(x,y)f(u)dy,
\end{equation}
where $G(x,y)$ is the Green's function of $(-\Delta)^k, n>2k$ corresponding to the Dirichlet problem. This integral equation is closely related to the higher order differential equation
\begin{equation*}
    \begin{cases}
        (-\Delta)^k u = f(u) &\text{ on } \halfrn,\\
        u=\frac{\partial u}{\partial x_n}=\cdots=\frac{\partial^{k-1}u}{\partial x_n^{k-1}}=0\;  &\text{ on } \partial \mathbb{R}^n_+.
    \end{cases}
\end{equation*}
Under suitable conditions for $f$, Lu and Zhu showed that the solution $u$ is axially symmetric with respect to some line parallel to the $x_n$-axis and $u(x)$ is non-decreasing in the $x_n$ direction.
Lu and Yang later \cite{LY} proved the following identity:

$$
  P_k u= x_n^{k + \frac{n}{2}} (- \Delta)^k \left( x_n^{k - \frac{n}{2}} u \right), x \in \mathbb{R}^n_+.
$$
With the help of the above identity, equation (\ref{P_k}) can also be reduced to an equivalent version on the upper half space. Hence Theorem \ref{thm1} implies the following axial symmetry result:

\begin{theorem}\label{cor1.4}
    Consider the equation
    \begin{equation}\label{Pk_half}
        (-\Delta)^kw = \frac{1}{|x_n|^{k+\frac{n}{2}}}f(|x_n|^{\frac{n}{2}-k}w)\; \text{ on } \halfrn.
    \end{equation}
    If $w\in C^{\infty}_0(\halfrn)$ is a positive solution to (\ref{Pk_half}), then
    there exists some parallel line to $x_n$-axis such that $w$ is symmetric with respect to such line.
\end{theorem}

\medskip

A natural question unsolved in \cite{LLY} is: does our moving plane method work for fractional order semilinear equations on hyperbolic spaces? In the present paper we will give an affirmitive answer. We first recall the definition of fractional order operators on $\mathbb{H}^n$, using the Helgason-Fourier transform. Let
$$e_{\lambda,\zeta}(x)=\left(\frac{\sqrt{1-|x|^2}}{|x-\zeta|}\right)^{n-1+i\lambda}, \; x\in\mathbb{B}^n, \lambda\in\mathbb{R}, \zeta\in\mathbb{S}^{n-1}.$$
The Helgason-Fourier transform of a function $f$ on $\mathbb{H}^n$ (ball model) is defined as
$$\hat{f}(\lambda,\zeta)=\int_{\bn}f(x)e^{-\lambda,\zeta}(x)dV,$$
as long as the integral exists (see Section 2 for more details). Then the fractional Laplacian on $\mathbb{H}^n$ can be defined as

$$\widehat{\left(-\Delta_{\mathbb{H}^n}\right)^\alpha}u(\lambda,\zeta)=\left(\frac{(n-1)^2+\lambda^2}{4}\right)^\alpha\hat{u}(\lambda,\zeta).$$
With the help of Helgason-Fourier transform, one can easily reduce a differential equation on $\mathbb{H}^n $ into an integral equation and we are particularly interested in the following integral equation:

\begin{equation}\label{eq1.14}
    u(x)=\int_{\hs}G_\alpha(x,y)f(u)dV_y,
\end{equation}
where $G_\alpha$ is the Green’s function of the operator $(-\Delta_{\mathbb{H}^n}-\frac{(n-1)^2}{4})^{\alpha/2}$. When $f$ satisfies certain integrability assumptions, we have the following:

\begin{theorem}\label{thm1.5}
    Let $f$ be Lipschitz continuous, non-decreasing and $f^\prime(u)\in L^{\frac{n}{\alpha}}(\hs)$, $n\geq 3$ and $0<\alpha<3$. If for $q > \frac{n}{n - \alpha}$, $u\in W^{q,2}(\hs)$  is a positive solution to (\ref{eq1.14}), then $u$ is radially symmetric and strictly decreasing in the radial variable.
\end{theorem}

Our moving plane method towards integral equations on hyperbolic spaces relies on the precise knowledge of the asymptotic behavior of the corresponding Green's functions as well as their monotonicity.
Green's function estimates and heat kernel estimates on Riemannian manifolds have been one of the central problems in geometric analysis (see e.g. \cite{D-M, Li-Yau, Strichartz}).
Unfortunately, those results on general noncompact manifolds do not seem to be sufficient to study our problems.
On the other hand, hyperbolic spaces possess symmetry structure and can be understood as real symmetric spaces of rank one.
The Helgason-Fourier analysis theory on $\hs$ provides much more precise information of Green's functions and heat kernels and can hence be applied to solve many analysis problems on $\hs$ (see \cite{A-J}).
In fact, the Helgason-Fourier analysis plays a key role in the establishment of higher order Hardy-Sobolev-Maz'ya's inequalities (see \cite{L-Y,LY}). Simultaneously in \cite{LLY2}, with the help of Helgason-Fourier analysis theory, the authors obtained a series of Green's function estimates and further use them to establish some Hardy-Adams inequalities involving fractional order operators.
We realize that those Green's function estimates in \cite{LLY2} instantly imply a Hardy-Littlewood-Sobolev type inequality on $\mathbb{H}^n$, which is the key ingredient to perform our moving plane argument. As far as we know, Theorem \ref{thm1.5} is the first moving plane result of fractional order equations on $\mathbb{H}^n$.

\medskip

More generally, let $\mathcal{G}(x,y)$ be the Green's function of the operator $$\prod_{j=0}^{l-1}\left(-\Delta_{\hs}-\frac{(n-1)^2}{4}+\zeta_j^2\right)^{s_j/2},$$ where $\zeta_j \geq 0$.
We consider the following integral equation:
\begin{equation}\label{prod_eq}
    u(x)=\int_{\bn}\mathcal{G}(x,y)f(u)dV_y.
\end{equation}
\begin{theorem}\label{thm1.6}
    Denote $S_l=\sum_{j=0}^{l-1}s_j$ and let $f$ be Lipschitz continuous, non-decreasing and $f^\prime(u)\in L^{\frac{n}{S_l}}(\hs)$. Assume $0\leq s_j<3$ if $\zeta_j>0$ and $S_l<n$. If $u(x)$ is a positive solution to (\ref{prod_eq}), then $u$ is radially symmetric and strictly decreasing in the radial variable.
\end{theorem}

\begin{remark}
    It is necessary to require $0\leq s_j<3$ and $S_l<n$ in order to invoke Lemma \ref{k_1 est} and Lemma \ref{prod est}.
\end{remark}
The paper is organized as follows. In Section 2, we present some preliminaries including the basics of hyperbolic spaces and some tools which will be used throughout the article;
Section 3 provides crucial estimates of Green's functions for the differential operators under consideration. Our main theorems will be proved in Section 4.

\section{Notations and Preliminaries}
We first present some preliminaries concerning hyperbolic space which are needed in the sequel. More information is referred to \cite{Ahlfors, Ratcliffe}.

\subsection{Models of hyperbolic spaces}\label{sec2.1}
The hyperbolic $n-$space $\hs$ $(n\geq 2)$ is a complete simply connected Riemannian manifold with constant sectional curvature $-1$.
There are several analytic models of hyperbolic spaces, all of which are equivalent.
Among them, we describe two models here.

\begin{itemize}
    \item \textit{The Half-space model}:
    It is given by $\mathbb{R}^{n-1}\times \mathbb{R}^+=\{(x_1,\cdots,x_{n-1},x_n):x_n>0\}$,
    equipped with the Riemannian metric
    $$ds^2=\frac{dx_1^2+\cdots+dx_n^2}{x_n^2}.$$

    \item \textit{The Poincar\'e ball model}:
    It is given by the open unit ball $\mathbb{B}^n=\{x=(x_1,\cdots,x_n):x_1^2+\cdots+x_n^2<1\}\in\mathbb{R}^n$ equipped with the Poincar\'e metric
    $$ds^2=\frac{4\left(dx_1^2+\cdots+dx_n^2\right)}{\left(1-|x|^2\right)^2}.$$
    The distance from $x\in\mathbb{B}^n$ to the origin is $\rho(x)=\log\frac{1+|x|}{1-|x|}$.
    The hyperbolic volume element is $dV=\left(\frac{2}{1-|x|^2}\right)^2dx$.
    The hyperbolic gradient is $\nabla_{\mathbb{H}^n}=\frac{1-|x|^2}{2}\nabla$ and the associated Laplace-Beltrami operator is given by
    $$\Delta_{\mathbb{H}^n}=\frac{1-|x|^2}{4}\left((1-|x|^2)\Delta+2(n-2)\sum_{i=1}^nx_i\frac{\partial}{\partial x_i}\right).$$
\end{itemize}

\subsection{Sobolev spaces on hyperbolic spaces}\label{Subsection2.2}
We define the Sovolev space $W^{k,2}$ on the Poincar\'e ball model.
For any open set $\Omega\in\hs$ and $u\in C^{\infty}(\Omega)$, the $W^{k,2}$ norm of $u$ is defined to be
\begin{equation*}
    \|u\|_{W^{k,2}(\Omega)}=\sum_{0\leq j\leq k}\int_\Omega|\left(-\Delta_{\mathbb{H}^n}\right)^{\frac{j}{2}}u|^2dV,
\end{equation*}
where
\begin{equation*}
    |\left(-\Delta_{\mathbb{H}^n}\right)^{\frac{j}{2}}u|^2=
    \begin{cases}
        |\left(-\Delta_{\mathbb{H}^n}\right)^{\frac{j}{2}}u|^2 & \text{if } j \text{ is even,}\\
        |\nabla_{\mathbb{H}^n}(-\Delta_{\mathbb{H}^n})^{\frac{j-1}{2}}u|^2 & \text{if } j \text{ is odd.}
    \end{cases}
\end{equation*}
The Sobolov space $W^{k,2}(\Omega)$ is then the closure of $C^\infty(\Omega)$ with respect to $\|\cdot\|_{W^{k,2}(\Omega)}$.
As usual, $W^{k,2}_{0}(\Omega)$ denotes the closure of $C^\infty_{0}(\Omega)$ in ${W^{k,2}(\Omega)}$.
Liu \cite{Liu} established the following sharp Sobolev inequalities:
\begin{equation*}
    \int_{\bn}(P_ku)udV\geq S_{n,k}\left(\int_{\bn}|u|^{\frac{2n}{n-2k}} dV\right)^{\frac{n-2k}{n}},\; u\in C^{\infty}_0(\bn), 1\leq k\leq \frac{n}{2},
\end{equation*}
where $S_{n,k}$ is the sharp constant of the classical k-th order Sobolev inequality on $\mathbb{R}^n$.

\subsection{Fractional Laplacian}Recall that
the fractional Laplacian in $\rn$ is a nonlocal operator which can be defiend as a singular integral operator:
\begin{equation*}
    (-\Delta)^{\frac{\alpha}{2}}u(x)=C_{n,\alpha}\textup{p.v.}\int_{\rn}\frac{u(x)-u(z)}{|x-z|^{n+\alpha}}dz,
\end{equation*}
where $\alpha$ is a real number between $0$ and $2$, p.v. stands for the Cauchy principle value and
$$C_{n,\alpha}=\frac{2^\alpha\Gamma(\frac{n+\alpha}{2})}{\pi^{\frac{n}{2}}\Gamma(-\frac{\alpha}{2})}.$$
Equivalently, $(-\Delta)^{\frac{\alpha}{2}}$ can be defined in terms of Fourier transform:
\begin{equation*}
    (-\Delta)^{\frac{\alpha}{2}}u(x)=\mathcal{F}^{-1}[|2\pi\xi|^\alpha \mathcal{F}{u}(\xi)](x),
\end{equation*}
where $\mathcal{F}$ is the Fourier transform.

\subsection{The Helgason-Fourier transform on hyperbolic spaces}
In this subsection, we will introduce the Helgason-Fourier transform on hyperbolic spaces and define the fractional Laplacian on hyperbolic spaces using the Helgason-Fourier analysis. For complete details, we refer to \cite{H-F1,H-F2}. Let
$$e_{\lambda,\zeta}(x)=\left(\frac{\sqrt{1-|x|^2}}{|x-\zeta|}\right)^{n-1+i\lambda}, \; x\in\mathbb{B}^n, \lambda\in\mathbb{R}, \zeta\in\mathbb{S}^{n-1}.$$
Then the Fourier transform of a function $f$ on $\mathbb{H}^n$ (ball model) is defined as
$$\hat{f}(\lambda,\zeta)=\int_{\bn}f(x)e_{-\lambda,\zeta}(x)dV,$$
provided the integral exists. Moreover, the following inversion formula holds for $f\in C^\infty_0(\bn)$:
$$f(x)=D_n\int_{-\infty}^\infty\int_{\mathbb{S}^{n-1}}\hat{f}(\lambda,\zeta)e_{\lambda,\zeta}(x)|\mathfrak{c}(\lambda)|^{-2}d\lambda d\sigma,$$
where $D_n=(2^{3-n}\pi|\mathbb{S}^{n-1}|)^{-1}$ and $\mathfrak{c}(\lambda)$ is the Harish-Chandra's $\mathfrak{c}$-function given by
$$\mathfrak{c}(\lambda)=\frac{2^{n-1-i\lambda}\Gamma(n/2)\Gamma(i\lambda)}{\Gamma(\frac{n-1+i\lambda}{2})\Gamma(\frac{1+i\lambda}{2})}.$$
There also holds the Plancherel formula:
$$\int_{\bn}|f(x)|^2dV=D_n\int_{-\infty}^\infty\int_{\mathbb{S}^{n-1}}\hat{f}(\lambda,\zeta)|\mathfrak{c}(\lambda)|^{-2}d\lambda d\sigma.$$
Since $e_{\lambda,\zeta}(x)$ is an eigenfunction of $-\Delta_{\hs}$ with eigenvalue $\frac{(n-1)^2+\lambda^2}{4}$, we have for $f\in C^\infty_0(\hs)$,
$$\widehat{\Delta_{\hs}f}(\lambda,\zeta)=-\frac{(n-1)^2+\lambda^2}{4}\hat{f}(\lambda,\zeta).$$
Therefore, we define the fractional Laplacian on hyperbolic spaces as follows:
$$\widehat{\left(-\Delta_{\hs}\right)^\alpha}u(\lambda,\zeta)=\left(\frac{(n-1)^2+\lambda^2}{4}\right)^\alpha\hat{u}(\lambda,\zeta).$$

\subsection{Hardy-Littlewood-Sobolev Inequality}
The HLS inequality on hyperbolic $\bn$ is equivalent to the HLS inequality on the hyperbolic upper half spaces,
which was first proved for half spaces by Beckner \cite{Beckner}, and then for Poincar\'e ball by Lu and Yang \cite{L-Y}.

\medskip

\textbf{Theorem H.} \textit{Let $0<\lambda<n$ and $p=\frac{2n}{2n-\lambda}$. Then for $f,g\in L^p(\bn)$,}
\begin{equation}\label{HLS}
    \left|\int_{\bn}\int_{\bn}\frac{f(x)g(y)}{(2\sinh (\frac{\rho(T_y(x)}{2}))^\lambda}dV_xdV_y\right|\leq C_{n,\lambda}\|f\|_p\|g\|_p,
\end{equation}
\textit{where}
\begin{equation*}
    C_{n,\lambda}=\pi^{\lambda/2}\frac{\Gamma(\frac{n}{2}-\frac{\lambda}{2})}{\Gamma(n-\frac{\lambda}{2})}\left(\frac{\Gamma(\frac{n}{2})}{\Gamma(n)}\right)^{-1+\frac{\lambda}{n}}
\end{equation*}
\textit{is the best constant for the classical Hardy-Littlewood-Sobolev constant on $\mathbb{R}^n$.
Furthermore, the constant $C_{n,\lambda}$ is sharp and there is no nonzero extremal function for the inequality (\ref{HLS}).}

\subsection{Foliations of hyperbolic spaces}
A foliation is an equivalence relation on a manifold, the equivalence classes being connected, injectively  submanifolds, all of the same dimension.
Let $\mathbb{R}^{n,1}=(\mathbb{R}^{n+1},\cdot)$, where $\cdot$ is Lorentzian inner product defined by $x\cdot y=-x_0y_0+x_1y_1+\cdots+x_ny_n$. The hypberboloid model
of $\hs$ is the submanifold $\{x\in\mathbb{R}^{n,1}: x\cdot x=-1, x_0>0\}$. A particular directional foliation can be obtained by choosing any $x_i$ direction, $i=1,\cdots,n$.
Without loss of generality, we may choose $x_1$ direction. Denote $\mathbb{R}^{n,1}=\mathbb{R}^{1,1}\times\mathbb{R}^{n-1}$, where $(x_1,x_0)\in \mathbb{R}^{1,1}$. We define
$A_t=\tilde{A}_t\otimes Id_{\mathbb{R}^{n-1}}$, where $\tilde{A}_t$ is the hypberbolic rotation in $\mathbb{R}^{1,1}$,
$$\tilde{A}_t=\begin{pmatrix}
    \cosh t & \sinh t \\
    \sinh t & \cosh t
\end{pmatrix}.$$
Let $U=\hs\cap\{x_1=0\}$ and $U_t=A_t(U)$, $\hs$ is then foliated by $U_t$ and $\hs=\bigcup_{t\in\mathbb{R}}U_t$.
The reflection $I$ is an isometry such that $I^2=Id$ and $I$ fixes a hypersurface, by $I(x_0,x_1,x_2,\cdots,x_n)=(x_0,-x_1,x_2,\cdots,x_n)$.
Moreover, define $I_t=A_t\circ I\circ A_{-t}$, then $U_t$ is fixed by $I_t$.

\subsection{Hypergeometric functions}
We denote
$$F(a,b,c,d)=\sum_{k=0}^\infty\frac{(a)_k(b)_k}{(c)_k}\frac{z^k}{k!},$$
where $c_k$ is not equal to non-positive integers, and $(a)_k$ is the (rising) Pochhammer symbol,
which is defined by
\begin{equation*}
    (a)_k=\begin{cases}
        0,  &k=0\\
        a(a+1)\cdots(a+k-1), &k>0.
    \end{cases}
\end{equation*}

\section{Green's functions estimates}

In what follows, $a\lesssim b$ will stand for $a\leq Cb$ for some positive constant $C$ and $a\sim b$ will stand for $C^{-1}b\leq a\leq Cb$.

\subsection{Green's function estimates of $P_k$}
For convenience, we introduce some notations. Let

\begin{equation*}
    \theta(x, y)=
    \begin{cases}
        (R-|x|^2)(R-|y|^2)  &\text{if } x, y\in B_R,\\
        0  &\text{if } x\not\in B_R \text{ or } y\not\in B_R.
    \end{cases}
\end{equation*}
and
\begin{align*}
        &H:(0,\infty)\times[0,\infty)\longrightarrow\mathbb{R}:\\
        &H(s, t)=s^{k-\frac{n}{2}}\int_0^{\frac{t}{s}}\frac{z^{k-1}}{(z+1)^{\frac{n}{2}}}dz.
\end{align*}

In \cite{Bo}, Boggio gave a representation formula for the Green's function of the operator $(-\Delta)^k$ corresponding to the Dirichlet problem on the unit ball $B_1 \subset \mathbb{R}^n$:

$$
  G(x,y) = C(n,k) |x-y|^{2k - n} \int_0^{\frac{(1 - |x|^2)(1 - |y|^2)}{|x - y|^2}} \frac{z^{k-1}}{(z + 1)^{n/2}} dz,
$$
where $C(n,k)$ is a constant depending only on the dimension $n$ and order $k$. We will first establish a Boggio type formula on the concentric ball $B_R \subset B_1$ where $B_1$ is equipped with the Poincar\'e metric.

\begin{lemma}
Consider the Dirichlet problem

\begin{equation}
    \begin{cases}
        P_k u = f(u) &\text{ on } B_R,\\
        \nabla_{\mathbb{H}^n}^\alpha u = 0,\; |\alpha|\leq k-1 &\text{ on }\partial B_R.
    \end{cases}
\end{equation}
Denote $G_{\mathbb{H}^n} (x,y)$ to be the Green's function of the operator $P_k$, then we have the following formula:

\begin{align*}
    G_{\mathbb{H}^n}(x, y)&=C(n,k) \frac{\left(1-|x|^2\right)^{\frac{n}{2}-k}\left(1-|y|^2\right)^{\frac{n}{2}-k}}{|x-y|^{n-2k}}\int_0^{\frac{(R^2-|x|^2)(R^2-|y|^2)}{R^2|x-y|^2}}\frac{z^{k-1}}{(z+1)^{\frac{n}{2}}}dz\\
    &=C(n,k,R)\left(1-|x|^2\right)^{\frac{n}{2}-k}\left(1-|y|^2\right)^{\frac{n}{2}-k}H\left(R^2|x-y|^2,\theta(x, y)\right).
\end{align*}

\end{lemma}

\begin{proof}
  Recall the following formula proved in \cite{Liu}:

  $$P_ku=\left(\frac{1-|x|^2}{2}\right)^{k+\frac{n}{2}}(-\Delta)^k\left((\frac{1-|x|^2}{2})^{k-\frac{n}{2}}u\right).$$
  Since from the definition of Green's function,

  \begin{align*}
     u(x) &= \int_{B_R} P_k u(y) G_{\mathbb{H}^n} (x,y) dV_y, \\
            &= \int_{B_R} \left(\frac{1-|y|^2}{2}\right)^{k + \frac{n}{2}} (-\Delta)^k \left( \left(\frac{1-|y|^2}{2}\right)^{k - \frac{n}{2}} u\right) G_{\mathbb{H}^n} (x,y) dV_y.
  \end{align*}
  By letting $\tilde{u}(x) = \left(\frac{1-|x|^2}{2}\right)^{k - \frac{n}{2}} u(x)$, we have

  $$
    \left(\frac{1-|x|^2}{2}\right)^{-k+\frac{n}{2}} \tilde{u}(x) = \int_{B_R} (- \Delta)^k \tilde{u}(y) \left(\frac{1-|y|^2}{2}\right)^{k + \frac{n}{2}} G_{\mathbb{H}^n} (x,y) \left(\frac{1-|y|^2}{2}\right)^{-n} dy.
  $$
  Thus if we denote $G_R$ as the Green's function of $(-\Delta)^k$ on $B_R$, we have

  $$
    G_R(x,y) = \left(\frac{1-|x|^2}{2}\right)^{k - \frac{n}{2}} \left(\frac{1-|y|^2}{2}\right)^{k - \frac{n}{2}} G_{\mathbb{H}^n} (x,y).
  $$
  Then we obtain our conclusion after proper rescaling.
\end{proof}

Without causing any confusion, we will write the Green's function $G_{\hs}$ for $P_k$ on hyperbolic space as $G$ for short.
We will first derive some pointwise inequalities for the Green's function of the GJMS operator $P_k$ on the ball $B_R$,
with respect to the Dirichlet boundary conditions.
\begin{lemma}\label{lemma2}
    For all $s, t>0$, we have
    $$H_s(s, t)<0,\quad H_t(s, t)>0,\quad H_{st}(s, t)<0.$$
\end{lemma}
\begin{proof}
    We first perform the change of variable $z=\frac{y}{s}$, so that
    \begin{align*}
        H(s, t)&=s^{k-\frac{n}{2}}\int_0^t\frac{\left(\frac{y}{s}\right)^{k-1}}{\left(\frac{y}{s}-1\right)^{n/2}}\frac{dy}{s}\\
        &=\int_0^t\frac{z^{k-1}}{(z+s)^{n/2}}dz.
    \end{align*}
    Therefore,
    \begin{equation*}
        H_t(s, t)=\frac{t^{k-1}}{(t+s)^{n/2}}>0,  \;\; H_{st}(s, t)=-\frac{n t^{k-1}}{2(t+s)^{n/2}+1}<0
    \end{equation*}
    and
    \begin{equation*}
        H_s(s, t)=-\frac{n}{2}\int_0^t\frac{z^{m-1}}{(z+s)^{n/2+1}}dz<0.
    \end{equation*}
\end{proof}

In the following, we put $T_\lambda:=\{x\in\mathbb{H}^n: x_1=\lambda\}$ and $\Sigma_\lambda:=\{x\in B_R: x_1<\lambda\}$, for all $\lambda\in[0,R]$.
For any $x\in\mathbb{H}^n$, let $\bar{x}$ denote the reflection of $x$ about $T_\lambda$.

\begin{lemma}\label{lemma3}
    Let $\lambda\in [0,1)$, then for every $x\in B\cap T_\lambda$ and $y\in\Sigma_\lambda$, we have
    \begin{align*}
        &G_{x_1}(x, y)<0\\
        &G_{ x_1}(x, y)+G_{x_1}(x,\bar{y})\leq 0.
    \end{align*}
    Moreover, the second inequality is strict if $\lambda>0$.
\end{lemma}
\begin{proof}
  For abbreviation, we put $d:=R^2|x-y|^2=R^2|x-\bar{y}|^2>0, \theta=\theta(x, y)>0$ and $\bar{\theta}=\theta(x,\bar{y})>0$.
  Then
  \begin{align*}
      G(x,y)&=C\left[-\left(\frac{n}{2}-k\right)2x_1\left(1-|x|^2\right)^{n/2-k-1}(1-|y|^{n/2-k})H(d,\theta)\right.\\
      &\quad+\left.\left(1-|x|^2\right)^{n/2-k}\left(1-|y|^2\right)^{n/2-k}H_{x_1}(d,\theta)\right]<0,
  \end{align*}
  since
  \begin{equation*}
      H_{x_1}(d,\theta)=-2H_t(d,\theta)\left(1-|y|^2\right)x_1+2H_s(d,\theta)(x_1-y_1)<0, \text{ for } x_1\geq 0, x_1>y_1.
  \end{equation*}
Moreover,
  \begin{align*}
      &G_{x_1}(x,y)+G_{x_1}(x,\bar{y})\\
      %&\leq C\left(\left(1-|x|^2\right)^{n/2-k}\left(1-|y|^2\right)^{n/2-k}H_{x_1}(d,\theta)+\left(1-|x|^2\right)^{n/2-k}\left(1-|\bar{y}|^2\right)^{n/2-k}H_{x_1}(d,\bar{\theta})\right)\\%
      &\leq C\left(1-|x|^2\right)^{n/2-k}\left(H_s(d,\theta)(x_1-y_1)+H_s(d,\bar{\theta})(x_1-\bar{y}_1)\right.\\
      &\quad\left.-\left[H_t(d,\theta)\left(1-|y|^2\right)+H_t(d,\bar{\theta})\left(1-|\bar{y}|^2\right)\right]x_1\right)\\
      &\leq C\left[H_s(d,\theta)-H_s(d,\bar{\theta})\right](x_1-y_1)<0,
  \end{align*}
  where we used Lemma \ref{lemma2} and the fact that $\bar{\theta}<\theta$.
\end{proof}

\begin{lemma}\label{max G}
    \begin{align*}
        &G(x, y)>\max\{G(x,\bar{y}),G(\bar{x},y)\}\\
        &G(x, y)+G(x,\bar{y})> |G(x,\bar{y})-G(\bar{x},y)|.
    \end{align*}
\end{lemma}

\begin{proof}
    Concerning the first inequality, it suffices to prove $G(x,y)>G(x,\bar{y})$ due to symmetry,
    since $I_t$ is isometry and hence $d(x,\bar{y})=d(\bar{x},y)$. Moreover,
    \begin{equation*}
        d(x,y)=d(\bar{x},\bar{y})<d(x,\bar{y})=d(\bar{x},y).
    \end{equation*}
    Besides, since $|\bar{x}|>|x|, |\bar{y}|>|y|$, we have that
    \begin{equation*}
       \theta(x, y)>\theta(\bar{x},\bar{y}).
    \end{equation*}
    Thus, we may conclude that
    \begin{align*}
        G(x,y)&=C\left(1-|x|^2\right)^{\frac{n}{2}-k}\left(1-|y|^2\right)^{\frac{n}{2}-k}H\left(R^2|x-y|^2,\theta(x, y)\right)\\
        &>C\left(1-|x|^2\right)^{\frac{n}{2}-k}\left(1-|\bar{y}|^2\right)^{\frac{n}{2}-k}H\left(R^2|x-\bar{y}|^2,\theta(x, \bar{y})\right)\\
        &=G(x,\bar{y}).
    \end{align*}
\end{proof}

We now extend $u$ by zero if it is outside of $B_R$ and we define
\begin{equation*}
    \tilde{f}(s)=
    \begin{cases}
        f(s)\quad &\text{if } s>0,\\
        0\quad &\text{if } s=0.
    \end{cases}
\end{equation*}
We then provide some crucial estimates for directional derivatives which are related to the Hopf boundary lemma for second order problems.

\begin{lemma}
    Let $0<\lambda<R$, and suppose that $u(x)\geq u(\bar{x})$ for all $x\in\Sigma_\lambda$. Then $\frac{\partial u}{\partial x_1}<0$
    on $T_\lambda\cap B_R$.
\end{lemma}
\begin{proof}
    For all $x\in T_l\cap B_R$ we have
    $$\frac{\partial u}{\partial x_1}(x)=\int_{B_R}G_{x_1}(x,y)f(u(y))dy=\int_{\Sigma_\lambda}[G_{x_1}(x,y)f(u(y))+G_{x_1}(x,\bar{y})\tilde{f}(u(\bar{y}))]dy$$
    Since $\tilde{f}$ is non-decreasing, we have $f(u(y))\geq \tilde{f}(u(\bar{y}))\geq 0$ for all $y\in\Sigma_\lambda$. Moreover, $f(u(y))\not\equiv 0$ in $\Sigma_\lambda$,
    since otherwise $f(u)\equiv 0$ in $B_R$. However, this would imply $P_ku\equiv 0$, which contradicts the positivity of $u$. As a result,
    there exists a nonempty open set $\mathcal{O}_\lambda\subset \Sigma_\lambda$ such that $f(u(y))>\tilde{f}(u(\bar{y}))$ or $\tilde{f}(u(\bar{y}))>0$
    for all $y\in\mathcal{O}_\lambda$.
    Thus, by Lemma (\ref{lemma3}), if $f(u(y))>\tilde{f}(u(\bar{y}))$,
    $$\frac{\partial u}{\partial x_1}(x)<\int_{\Sigma_\lambda}(G_{x_1}(x,y)+G_{x_1}(x,\bar{y}))\tilde{f}(u(\bar{y}))dy\leq 0, \text{ for all } x\in T_\lambda\cap B_R,$$
    and if $\tilde{f}(u(\bar{y}))>0$,
    $$\frac{\partial u}{\partial x_1}(x)\leq\int_{\Sigma_\lambda}(G_{x_1}(x,y)+G_{x_1}(x,\bar{y}))\tilde{f}(u(\bar{y}))dy<0 \text{ for all } x\in T_\lambda\cap B_R.$$
    In any case, we have
    $$\frac{\partial u}{\partial x_1}(x)<0 \text{ for all } x\in T_\lambda\cap B_R.$$
\end{proof}

Now we show this estimate is still true if we move $T_\lambda$ to the origin a little bit farther.

\begin{lemma}\label{dir deriv}
    Let $0<\lambda<R$, and suppose that $u(x)\geq u(\bar{x})$ for all $x\in\Sigma_\lambda$. Then there exists $\gamma\in(0,\lambda)$ such that $\frac{\partial u}{\partial x_1}<0$
    on $T_l\cap B_R$ for all $l\in(\lambda-\gamma,\lambda)$.
\end{lemma}

Before proving the above lemma, we need the following result which may be viewed as the higher order analogue to the Hopf lemma in hyperbolic spaces.

\begin{lemma}\label{higher order hopf}
    If $x_0\in\partial B_R$ and $\nu$ is a unit vector with $\nu\cdot x_0<0$, then $\frac{\partial^k u}{\partial\nu^k}(x_0)>0$.
\end{lemma}
\begin{proof}
   Let $$v(x)=\left(\frac{1-|x|^2}{2}\right)^{k-\frac{n}{2}}u(x).$$
   It is known that $\frac{\partial^kv}{\partial \nu^k}(x_0)>0$ and
   $$\frac{\partial^k v}{\partial \nu^k}=\sum_{i=0}^k\binom{k}{i}\frac{\partial^i}{\partial \nu^i}\left(\frac{1-|x|^2}{2}\right)^{k-\frac{n}{2}}\frac{\partial^{k-i}}{\partial \nu^{k-i}}u(x).$$
   We also have $(\frac{\partial}{\partial \nu})^m u(x_0)=0$ for all $m=0,\cdots,k-1$, from the boundary conditions. Thus,
   $$\left(\frac{\partial }{\partial \nu}\right)^kv(x_0)=\left(\frac{1-|x|^2}{2}\right)^{k-\frac{n}{2}}\left(-\frac{\partial}{\partial \nu}\right)^ku(x)>0,$$
   then $\left(\frac{\partial}{\partial \nu}\right)^ku(x)>0$.
\end{proof}

\begin{proof}[Proof of Lemma \ref{dir deriv}]
    In view of Lemma \ref{higher order hopf}, we know that for any $x_0\in T_\lambda\cap\partial B_R$,
    $$(-1)^k\left(\frac{\partial}{\partial x_1}\right)^{k-1}\frac{\partial u}{\partial x_1}(x_0)=\left(-\frac{\partial}{\partial x_1}\right)^{k}u(x_0)>0.$$
    We also know that $(\frac{\partial}{\partial x_1})^m u(x_0)=0$ for all $m=0,\cdots,k-1$, from the boundary conditions. Thus, there exists $a=a(x_0)>0$ such that
    \begin{equation}\label{eq8}
        \frac{\partial u}{\partial x_1}(x)<0, \text{ for all } x\in\mathcal{U}_a(x_0)\cap B_R,
    \end{equation}
    where $\mathcal{U}_a(y):=\{x\in\hs:\max\limits_{1\leq i\leq n} \rho(x_i,y_i)<a\}$.
    Then by the compactness of $T_\lambda\cap\partial B_R$, there exsits an uniform $\bar{a}>0$ such that
    \begin{equation*}
        \frac{\partial u}{\partial x_1}(x)<0, \text{ for all } x\in A:=\bigcup_{x_0\in T_\lambda\cap\,\partial B_R}\mathcal{U}_{\bar{a}}(x_0)\cap B_R.
    \end{equation*}
    Now let $ K:= (T_\lambda \cap B)\setminus A$ and consider $K_d:=K-de_1$ for $d>0$, there exists $\delta>0$ such that $\frac{\partial u}{\partial x_1}<0$ on $K_d$
    for all $d\in[0,\delta]$, due to Lemma \ref{higher order hopf}. Finally, let $\gamma:=\min\{\bar{a},\delta\}$, the statement of Lemma \ref{dir deriv} follows.
\end{proof}

\subsection{Estimates of $k_\alpha$}
We denote $k_\alpha=\left(-\Delta_{\hs}-\frac{(n-1)^2}{4}\right)^{-\alpha/2}$. It is known in \cite{LLY2} that $k_\alpha$ has the following estimate.
\begin{lemma}
Let $n\geq 3, 0<\alpha<3$. There holds
\begin{equation}\label{k_est}
    k_\alpha(\rho)=\frac{1}{\gamma(\alpha)}\cdot\frac{1}{\rho^{n-\alpha}}+O\left(\frac{1}{\rho^{n-\alpha-1}}\right), \quad 0<\rho<1,
\end{equation}
where $$\gamma(\alpha)=\frac{\pi^{\frac{n}{2}}2^\alpha\Gamma\left(\frac{\alpha}{2}\right)}{\Gamma\left(\frac{n}{2}-\frac{\alpha}{2}\right)}.$$
\end{lemma}

\subsection{Estimates of $\mathcal{G}(x,y)$}
We start by recalling that $\mathcal{G}(x,y)$ is the Green’s function of operator $\prod_{j=0}^{l-1}\left(-\Delta_{\hs}-\frac{(n-1)^2}{4}+\zeta_j^2\right)^{s_j/2}$.
For simplicity, we denote each term in the product by
\begin{equation*}
    k_{\zeta_j,s_j}=\left(-\Delta_{\hs}-\frac{(n-1)^2}{4}+\zeta_j^2\right)^{-s_j/2}.
\end{equation*}
Then the following estimate has been proved by Li, Lu and Yang \cite{LLY2}:
\begin{lemma}\label{k_1 est}
    Let $n\geq 3$ and $0<s_j<n$. There exists $0<\epsilon_1<\min\{1,n-s_j\}$ such that
    \begin{equation*}
        k_{\zeta_j,s_j}\leq \frac{1}{\gamma(s_j)}\cdot\frac{1}{\rho^{n-s_j}}+O\left(\frac{1}{\rho^{n-s_j-\epsilon_1}}\right),\; 0<\rho<1, \zeta_j>0.
    \end{equation*}
\end{lemma}
We will also need the following asymptotic estimate (see \cite{L-Y} Lemma 2.4).
\begin{lemma}\label{conv}
    Let $0<\alpha<n, 0<\beta<n,0<\alpha+\beta<n$, and $\lambda_1+\lambda_2>\alpha+\beta-2$. Then for $0<\rho<1$,
   \begin{align*}
    &\frac{1}{(2\sinh\frac{\rho}{2})^{n-\alpha}(\cosh\frac{\rho}{2})^{\lambda_1}}* \frac{1}{(2\sinh\frac{\rho}{2})^{n-\beta}(\cosh\frac{\rho}{2})^{\lambda_2}}\\
    &\leq \frac{\gamma(\alpha)\gamma(\beta)}{\gamma(\alpha+\beta)}\frac{1}{\rho^{n-\alpha-\beta}}+O\left(\frac{1}{\rho^{n-\alpha-\beta-\epsilon}}\right),
   \end{align*}
   where $\epsilon$ satisfies $0<\epsilon\leq 1$ if $0<\alpha+\beta<n-1$ and $0<\epsilon<n-\alpha-\beta$ if $n-1\leq \alpha+\beta<n$.
\end{lemma}
Now we give the estimate for the product of $k_{\zeta_j,s_j}$'s, which may also be of independent interest. Recall $S_l=\sum_{j=0}^{l-1}s_j$, we get the following,
\begin{lemma}\label{prod est} Let $n\geq 3$. Assume each $\zeta_j\geq 0$ and each $s_j\in (0,3)$ if $\zeta_j>0, S_l<n$, there exists $\epsilon_2\in(0,\min\{1,n-S_l\})$ such that
    \begin{align}\label{pr_est}
        \prod_{j=0}^{l-1}\left(-\Delta_{\hs}-\frac{(n-1)^2}{4}+\zeta_j^2\right)^{-s_j/2}
        \leq \frac{1}{\gamma(S_l)}\cdot\frac{1}{\rho^{n-S_l}}+O\left(\frac{1}{\rho^{n-S_l-\epsilon_1}}\right),\; 0<\rho<1.
    \end{align}
\end{lemma}
\begin{proof}
We prove it by induction. (\ref{pr_est}) is vaild for $k=1$ due to Lemma \ref{k_1 est}. Assume it is valid for $l>1$, then for $l+1$, we have
\begin{align*}
    &\prod_{j=0}^{l}\left(-\Delta_{\hs}-\frac{(n-1)^2}{4}+\zeta_j^2\right)^{-s_j/2}
    =\prod_{j=0}^{l-1}\left(-\Delta_{\hs}-\frac{(n-1)^2}{4}+\zeta_j^2\right)^{-s_j/2}*l_{\zeta_k,s_l}\\
    &\leq \left\{\frac{1}{\gamma(S_{k+1})}\cdot\frac{1}{\rho^{n-S_{l+1}}}+O\left(\frac{1}{\rho^{n-S_{l+1}-\epsilon_1}}\right)\right\}*k_{\zeta_k,s_k}\\
    &\lesssim \left\{\frac{1}{\gamma(S_k)}\frac{\left(\cosh\frac{\rho}{2}\right)^{-2\hat{\zeta}-S_l+1}}{\left(2\sinh(\frac{\rho}{2})\right)^{n-S_l}}+O\left(\frac{(\cosh\frac{\rho}{2})^{-2\hat{\zeta}-S_l+2}}{(\sinh\frac{\rho}{2})^{n-S_l+1}}\right)\right\}*\\
    &\quad\; \left\{\frac{1}{\gamma(s_l)}\frac{\left(\cosh\frac{\rho}{2}\right)^{-2\zeta_l^\prime-s_l+1}}{\left(2\sinh(\frac{\rho}{2})\right)^{n-s_l}}+O\left(\frac{(\cosh\frac{\rho}{2})^{-2\zeta^\prime_l-s_l+2}}{(\sinh\frac{\rho}{2})^{n-s_l+1}}\right)\right\}\\
    &\lesssim\frac{1}{\gamma(S_{l+1})}\cdot\frac{1}{\rho^{n-S_{l+1}}}+O\left(\frac{1}{\rho^{n-S_{l+1}-\epsilon_1}}\right).
\end{align*}
where $0<\hat{\zeta}<\sum_{j=0}^{l-1}\zeta_j$. If we choose $\epsilon<\hat{\zeta}$, then by Lemma \ref{conv}, there exists $0<\epsilon_1<\min\{1,n-S_k\}$ such that
(\ref{pr_est}) is true for $l+1$.
\end{proof}

\subsection{Monotonicity of Green’s functions}
Our final goal in this section is to prove that Green's functions of the aforementioned operators are all decreasing with respect to the geodesic distance.
\\~\\
We first recall that the heat kernel on $\bn$, denoted by $e^{t\Delta_{\hs}}$, is given explicitly by the following formulae (see, e.g. \cite{D-M, G-N}):
\begin{itemize}
    \item If $n=2m+1$, then
        \begin{equation}\label{heat1}
            e^{t\Delta_{\hs}}=2^{-m-1}\pi^{-m-\frac{1}{2}}t^{-\frac{1}{2}}e^{-\frac{(n-1)^2}{4}t}\left(-\frac{1}{\sinh\rho}\frac{\partial}{\partial\rho}\right)^me^{-\frac{\rho^2}{4t}}.
        \end{equation}
    \item  If $n=2m$, then
        \begin{equation}\label{heat2}
            e^{t\Delta_{\hs}}=(2\pi)^{-m-\frac{1}{2}}t^{-\frac{1}{2}}e^{-\frac{(n-1)^2}{4}t}\int_\rho^{\infty}\frac{\sinh r}{\sqrt{\cosh r-\cosh\rho}}\left(-\frac{1}{\sinh r}\frac{\partial}{\partial r}\right)^me^{-\frac{r^2}{4t}}dr.
        \end{equation}
\end{itemize}

\begin{lemma}
    Let $G_\alpha$ be the Green's function of the operator $\left(-\Delta_{\hs}-\frac{(n-1)^2}{4}\right)^{\alpha/2}$ on the hypberbolic ball $\bn$.
    Then for any fixed $y\in\bn$, $G_\alpha(x,y)$ is a positive radially decreasing function with respect to the geodesic distance $\rho=d(x,y)$.
\end{lemma}

\begin{proof}\label{ka}
    By the Mellin type expression, for $0<\alpha<\min\{3,n\}$,
    \begin{equation}\label{M-type}
        k_\alpha(\rho)=\left(-\Delta_{\hs}-\frac{(n-1)^2}{4}\right)^{-\alpha/2}=\frac{1}{\Gamma(\alpha/2)}\int_0^\infty t^{\frac{\alpha}{2}-1}e^{t\Delta_{\hs}+(n-1)^2/4}dt.
    \end{equation}
    If $n=2m+1$, by (\ref{heat1}), we have
    \begin{align*}
        &k_{\alpha,2m+1}(\rho):=k_\alpha(\rho)=C_1\int_0^\infty t^{\frac{\alpha-3}{2}}\left(-\frac{1}{\sinh\rho}\frac{\partial}{\partial \rho}\right)^me^{-\frac{\rho^2}{4t}}dt,
    \end{align*}
    where $C_1=\frac{2^{-m-1}\pi^{-m-1/2}}{\Gamma(\alpha/2)}$ is constant. Thus,
    \begin{align*}
        &\frac{d}{d\rho}k_{\alpha,2m+1}(\rho)=C_1(-\sinh\rho)\int_0^\infty t^{\frac{\alpha-3}{2}}\left(-\frac{1}{\sinh\rho}\frac{\partial}{\partial\rho}\right)^{m+1}e^{-\frac{\rho^2}{4t}}dt\\
        &=C_1(-\sinh\rho)\cdot k^{2m+3}_\alpha(\rho)<0.
    \end{align*}
    If $n=2m$, by (\ref{heat2}) and Fubini's theorem, we have
    \begin{align*}
        k_{\alpha,2m}(\rho):=k_\alpha(\rho)&=\sqrt{2}\int_\rho^{\infty}\frac{\sinh r}{\sqrt{\cosh r-\cosh\rho}}\left(\int_0^\infty t^{\frac{\alpha-3}{2}}\left(-\frac{1}{\sinh r}\frac{\partial}{\partial r}\right)^me^{-\frac{r^2}{4t}}dt\right)dr\\
        &=\sqrt{2}\int_\rho^\infty\frac{\sinh r}{\sqrt{\cosh r-\cosh\rho}}\cdot k_{\alpha,2m+1}(r)dr.
    \end{align*}
    Now if we set $s=\sqrt{\cosh r-\cosh\rho}=\sqrt{2\sinh^2(\frac{r}{2})-2\sinh^2(\frac{\rho}{2})}$, we get
    \begin{align*}
        &\frac{d}{d\rho}k_{\alpha,2m}(r)=\frac{\sqrt{2}}{2}\int_0^\infty \frac{d}{d\rho}k_{\alpha,2m+1}(r(s,\rho))ds<0
    \end{align*}
    becasue
    \begin{align*}
        \frac{d}{d\rho}k_{\alpha,2m+1}=\frac{dk_{\alpha,2m+1}(r)}{dr}\frac{dr}{d\rho}=\frac{dk_{\alpha,2m+1}(r)}{dr}\frac{\sinh\rho}{\sinh r}<0.
    \end{align*}
\end{proof}

\begin{lemma}
    Let $G_{\alpha,\zeta}$ be the Green's function of the operator $\left(-\Delta_{\hs}-\frac{(n-1)^2}{4}+\zeta^2\right)^{\alpha/2}$ on the hypberbolic ball $\bn$.
    Then for any fixed $y\in\bn$, $G_{\alpha,\zeta}(x,y)$ is a positive radially decreasing function with respect to the geodesic distance $\rho=d(x,y)$.
\end{lemma}
\begin{proof}
In order to show that $k_{\alpha,\zeta}(\rho)$ is also decreasing, we only need to notice that $k_{\alpha,\zeta}(\rho)$ is  defined as follows,
\begin{align*}
    &k_{\alpha,\zeta,2m+1}(\rho):=k_{\alpha,\zeta}(\rho) = C \int_0^\infty t^{\frac{\alpha}{2}-1}\left[\left(-\frac{1}{\sinh\rho}\frac{\partial}{\partial\rho}\right)^me^{-\frac{\rho^2}{4t}}\right]e^{-t\zeta^2}dt, \text{ if } n=2m+1,\\
    &k_{\alpha,\zeta,2m}(\rho):=k_{\alpha,\zeta}(\rho)\\
    &= C \int_\rho^{\infty}\frac{\sinh r}{\sqrt{\cosh r-\cosh\rho}}\left[\left(\int_0^\infty t^{\frac{\alpha-3}{2}}\left(-\frac{1}{\sinh r}\frac{\partial}{\partial r}\right)^me^{-\frac{r^2}{4t}}\right]e^{-t\zeta^2}dt\right)dr, \text{ if } n=2m.\\
\end{align*}
It is obvious to see that the derivative of $k_{\alpha,\zeta}$ with respect to $\rho$ is independent of $\zeta$ and the rest follows directly the proof of Lemma \ref{ka} so we omit it.
\end{proof}

\begin{lemma}
    Let $\mathcal{G}(x,y)$ be the Green's function of the operator  $\prod_{j=0}^{l-1}\left(-\Delta_{\hs}-\frac{(n-1)^2}{4}+\zeta_j^2\right)^{s_j/2}$ on the hypberbolic ball $\bn$.
    $S_l<n$ with $0\leq s_j<3$. Then for any fixed $y\in\bn$, $\mathcal{G}(x,y)$ is a positive radially decreasing function with respect to the geodesic distance $\rho=d(x,y)$.
\end{lemma}
\begin{proof}
    It suffices to show that the convolution of two positive radially decreasing functions $H_1, H_2$, is still positve radially decreasing (see also \cite{LLY}).
    Denote $$L(x,y)=\int_{\bn}H_(x,z)H_2(z,y)dV_z.$$
    We first show $L(x,y)=L(\rho)$. Since for any isometry $T: \bn\rightarrow\bn$,
    \begin{align*}
        L(Tx,Ty)&=\int_{\bn}H_1(Tx,z)H_2(z,Ty)dV_z\\
        &=\int_{\bn}H_1(x,T^{-1}z)H_2(T^{-1}z,y)dV_z\\
        &=\int_{\bn}H_1(x,z)H_2(z,y)dV_z\\
        &=L(x,y).
    \end{align*}
    Now for fixed $y$, consider geodesic ray from $y$ and two points $x$ and $\bar{x}$ on the ray. Without loss of generality, we may assume $\rho(x,y)\leq\rho(\bar{x},y)$. Due to the foliation of $\bn$,
    there exists an unique reflection $I_t$ such that $I_t(x)=I_t(\bar{x})$. Then it is easy to show that $L(x,y)-L(\bar{x},y)\geq 0$ and the lemma follows.
\end{proof}

\section{Hyperbolic symmetry of the solutions}
\subsection{Proof of Theorem \ref{thm1}}

We will first neeed to show that the differential equation (\ref{P_k}) is equivalent to an integral equation using the Helgason-Fourier trasnformation on hypberbolic spaces.

\begin{lemma}\label{lemma1}
    If $u\in W^{k,2}_0(\hs)$ is a positive weak solution of the higher order
    differential equation (\ref{P_k}), then $u$ must satisfy the following integral equation for the
    Green's function $G(x, y)$ of the operator $P_k$ on the hyperbolic ball $\mathbb{B}^n$:
    \begin{equation}\label{inteq1}
        u(x)=\int_{\mathbb{B}^n} G(x,y)f(u)dV_y.
    \end{equation}
\end{lemma}

\begin{proof}
    $u\in W^{k,2}_0(\hs)$ is a weak solution of \eqref{inteq1} if and only if for any $\phi\in C^\infty_0(\hs)$,
    $$\inth \sum_{j=1}^kc_j\nabla^j_{\hs}u\nabla^j_{\hs}\phi = \inth f(u)\phi dV.$$
    Using Helgason-Fourier transform on $\hs$, such definition of weak solution is equivalent to
    \begin{equation*}
        D_n\int_{-\infty}^{+\infty}\int_{\mathbb{S}^{n-1}}\left[\prod_{j=1}^k\left(\frac{\tau^2+(2k-1)^2}{4}\right)\right]\hat{u}(\tau,\sigma)\hat{\phi}(\tau,\sigma)|\mathfrak{c}(\tau)|^{-2}d\sigma d\tau = \inth f(u)\phi dV.
    \end{equation*}
    Let $\psi$ satisfy $P_k\psi = \phi$, that is, $\psi(x)=\inth G(x, y)\phi(y)dV_y$. Under Helgason-Fourier transform, we have
    \begin{equation*}
        \hat{\phi}(\tau,\sigma)=\prod_{j=1}^k\frac{\tau^2+(2k-1)^2}{4}\hat{\psi}(\tau,\sigma).
    \end{equation*}
    Now if we replace $\phi$ by $\psi$, we get
    \begin{equation*}
        D_n\int_{-\infty}^{+\infty}\int_{\mathbb{S}^{n-1}}\hat{u}(\tau,\sigma)\hat{\psi}(\tau,\sigma)|c(\tau)|^{-2}d\sigma d\tau = \inth f(u)\left(\inth G(x,y)\phi(y)dV_y\right)dV_x.
    \end{equation*}
    Applying Plancherel formula to the left hand side, we get
    \begin{equation*}
        \inth u(x)\phi(x)dV = \inth\left(\inth G(x,y)f(x)dV_x\right)\phi(y)dV_y,
    \end{equation*}
    which is true for any $\phi\in\cs(\hs)$. This immediately implies that a solution of the differential equation (\ref{P_k})
    is a solution of the integral equation (\ref{inteq1}). The other direction of the statement can be shown in the same spirit and we omit it.
    \end{proof}

    In particular, it has been shown in \cite{LY}, that the Green's functions of $P_k$ satisfies
    \begin{equation*}
    P_k^{-1}(\rho)=\frac{\Gamma(\frac{n}{2})}{2^n\pi^{\frac{n}{2}}\Gamma(k)\Gamma(k+1)}\frac{\left(\cosh\frac{\rho}{2}\right)^{-n}}{\left(\sinh\frac{\rho}{2}\right)^{n-2k}}F\left(k-\frac{n-2}{2},k;k+1;\cosh^{-2}\frac{\rho}{2}\right),
    \end{equation*}
    where $F$ is the hypergeomtric function.
    Furthermore, we have, for $1\leq k < \frac{n}{2}$,
    \begin{equation}\label{pk_est}
    P_k^{-1}(\rho)\leq\frac{1}{\gamma_n(2k)}\left[\left(\frac{1}{2\sinh\frac{\rho}{2}}\right)^{n-2k}-\left(\frac{1}{2\cosh\frac{\rho}{2}}\right)^{n-2k}\right], \;\rho>0,
    \end{equation}

    where $\rho=\log\frac{1+|x|}{1-|x|}$ is the hyperbolic distance from $x\in \mathbb{B}^n$ to the origin.
    As a result, we have the following proposition immediately.
    \begin{prop}
    Let $G(x,y)$ be the Green's function of the operator $P_k$ on the hyperbolic ball $\mathbb{B}^n$.
    Then for any fixed $y\in \mathbb{B}^n, G(x,y)$ is a positive radially decreasing function with respect
    to the geodesic distance $\rho=d(x,y)$.
\end{prop}
    Due to Lemma \ref{lemma1}, it suffices to prove the symmetry of the solution to the integral equation (\ref{inteq1}).
    Recall $U=\hs\cap\{x_1=0\}$ and  $U_\lambda=A_\lambda(U), \Sigma_\lambda=\cup_{s<\lambda}U_s$.
    For any $x\in\Sigma_\lambda$, denote $\bar{x}=I_\lambda(x)$ and $u_\lambda(x)=u(\bar{x})$. Consider $\lambda>0$,
    \begin{align*}
    &u(x)-u_\lambda(x)\\
    &=\int_{\Sigma_\lambda} G(x,y)f(u)dV_y+\int_{\Sigma^c_\lambda} G(x,y)f(u)dV_y\\
    &\quad -\int_{\Sigma_\lambda} G(x,y)f(u_\lambda)dV_y-\int_{\Sigma^c_\lambda} G(x,y)f(u_\lambda)dV_y\\
    &=\int_{\Sigma_\lambda} G(x,y)f(u)dV_y+\int_{\Sigma_\lambda} G(x,\bar{y})f(u_\lambda)dV_y\\
    &\quad -\int_{\Sigma_\lambda} G(x,y)f(u_\lambda)dV_y-\int_{\Sigma_\lambda} G(x,\bar{y})f(u_\lambda)dV_y\\
    &=\int_{\Sigma_\lambda} G(x,y)f(u)dV_y+\int_{\Sigma_\lambda} G(\bar{x},y)f(u_\lambda)dV_y\\
    &\quad -\int_{\Sigma_\lambda} G(x,y)f(u_\lambda)dV_y-\int_{\Sigma_\lambda} G(\bar{x},y)f(u)dV_y\\
    &=\int_{\Sigma_\lambda} \left(G(x,y)-G(\bar{x},y)\right)\left(f(u)-f(u_\lambda)\right)dV_y
    \end{align*}
    In $\Sigma_\lambda$, we denote $\Sigma_\lambda^-=\{x\in\Sigma_\lambda: u_\lambda(x)>u(x)\}$, and we would like to show that $\Sigma_\lambda^-$ is of measure zero.
    In fact, we have from (\ref{pk_est}) that,
    \begin{align*}
    u_\lambda(x)-u(x)&\leq \int_{\Sigma_\lambda^-}G(x,y)\left(f(u_\lambda)-f(u)\right)dV_y\\
    &\leq C\int_{\Sigma_\lambda^-}\left(\frac{1}{2\sinh \frac{d(x,y)}{2}}\right)^{n-2k}\left(f(u_\lambda)-f(u)\right)dV_y.
    \end{align*}
    Moreover, by the Hardy-Littlewood-Sobolev inequality on $\hs$, mean value theorem and H\"older's inequality, we have
    \begin{align*}
    \|u-u_\lambda\|_{L^q(\Sigma_\lambda)}&\leq C\|f(u)-f(u_\lambda)\|_{\frac{nq}{n+2kq}}\\
    &=C\left(\int_{\Sigma_\lambda^-}\left(f^\prime(\xi)(u(y)-u_\lambda(y))\right)^{\frac{nq}{n+2kq}}\right)^{\frac{1}{q}+\frac{2k}{n}}\\
    %&\leq C\left(\int_{\Sigma_t^-}f^\prime(\xi)^\frac{n}{2k}\right)^{\frac{2k}{n}}\left(\int_{\Sigma_t^-}\left(u(y)-u_t(y)\right)^q\right)^\frac{1}{q}\\%
    &\leq C\|f^\prime(\xi)\|_{L^{\frac{n}{2k}}(\Sigma_\lambda^-)}\|u-u_\lambda\|_{L^q(\Sigma_\lambda^-)},
    \end{align*}
    for some $\frac{n}{n-2k}<q<\infty$.
    Here $\xi=\theta u + (1-\theta)u_\lambda$ with $0<\theta<1$.
    By monotonicity of $f^{\prime}$, we are able to pick $\lambda$ large enough so that $C\|f^\prime(\xi)\|_{\frac{n}{2k}}<1$.
    This implies that for large enough $\lambda, \|u-u_\lambda\|_{L^q(\Sigma_\lambda^-)}=0$, which shows that $\Sigma_\lambda^-$ is of measure zero.
    \\~\\
    Now we shift $U_\lambda$ as long as $u\geq u_\lambda$ in $\Sigma_\lambda$. Suppose that there exists such $\bar{\lambda}$ that $u(x) > u_{\bar{\lambda}}(x)$ on $\Sigma_{\bar{\lambda}}$.
    We deduce again by a standard compactness argument that there exist $\bar{\lambda}-\varepsilon<\lambda\leq \bar{\lambda}$ such that
    $$\|u-u_\lambda\|_{L^q(\Sigma_\lambda)}\leq C\|f^\prime(\xi)\|_{L^{\frac{n}{x2k}}(\Sigma_\lambda^-)}\|u-u_\lambda\|_{L^q(\Sigma_\lambda^-)}.$$
    When $\varepsilon$ is small, $\Sigma_\lambda^-$ is close to zero so that $\|u_t-u\|_{L^q(\Sigma_\lambda^-)}=0$.
    This imples we can keep moving the plane $U_{\lambda_0}$. Now we see $u(x)\leq u(\bar{x})$ with respect to $\Sigma_{\lambda_0}$.
    A similar argument shows $u(\bar{x})\geq u(x)$ by rotation. Consequently, there exsits a point $p$ such that $u$ is constant on the geodesic spheres center at $p$.

\subsection{Proof of Theroem \ref{thm2}}

Now we start the moving plane procedure by shifting the plane $T_\lambda$ from the initial tangent position $T_R$ towards the interior of $B_R$.

\begin{lemma}\label{lemma 3.8}
    There exists $\varepsilon>0$ such that for all $\lambda\in[R-\varepsilon,R)$ we have
    \begin{equation}\label{eq9}
    u(x)>u(\bar{x}) \text{ for } x\in\Sigma_\lambda, \;\;
    \frac{\partial u}{\partial x_1}(x)<0 \;\text{ for } x\in T_\lambda\cap B_R.
    \end{equation}
\end{lemma}
\begin{proof}
    Since $T_R\cap\,\partial B_R=\{e_1\}$, where $e_1=(R,0,\cdots,0)$. We conclue from (\ref{eq8}) that there exists $\varepsilon>0$ such that $\frac{\partial u}{\partial x_1}(x)<0$, for
    $x\in B_R\setminus\Sigma_{R-2\varepsilon}$. As a result, (\ref{eq9}) holds for all $\lambda\in[R-\varepsilon, R]$.
\end{proof}

We then finish the proof of the theorem by moving the plane until it reaches the origin.

Let $$\Lambda:=\{\lambda\in (0, R): u(x)>u(\bar{x}), \forall x\in\Sigma_\lambda,\frac{\partial u}{\partial x_1}(x)<0, \forall x\in T_\lambda\cap B\}.$$
By Lemma \ref{lemma 3.8}, we know that
$[R-\varepsilon, R)\subset \Lambda$. Let $\bar{\lambda}\in [0, R)$ be the smallest number such that $(\bar{\lambda}, R)\subset\Lambda$. Once agian, we would like to  show $\bar{\lambda}=0$
so that $\Lambda=(0, R)$. Assume by contradiction that $\bar{\lambda}>0$. Then according to Lemma $\ref{dir deriv}$ we get that
\begin{equation}\label{cty}
\text{there exists } \gamma\in(0,\bar{\lambda}) \text{ such that } \frac{\partial u}{\partial x_1}<0, \text{ on } T_l\cap B_R,\; \forall l\in(\bar{\lambda}-\gamma,\bar{\lambda}).
\end{equation}
Now for all $x\in\Sigma_{\bar{\lambda}}$,
\begin{align*}
    &u(x)-u(\bar{x})=\int_{B_R}(G(x,y)-G(\bar{x},y))f(u(y))dy\\
    &=\int_{\Sigma_{\bar{\lambda}}} G(x,y)f(u(y))dy+\int_{\Sigma^c_{\bar{\lambda}}} G(x,y)f(u(y))dy-\int_{\Sigma_{\bar{\lambda}}} G(x,y)f(u(y)))dy-\int_{\Sigma^c_{\bar{\lambda}}} G(x,y)f(u(y))dy\\
    &=\int_{\Sigma_{\bar{\lambda}}}(G(x,y)-G(\bar{x},y))f(u(y))dy+\int_{\Sigma_{\bar{\lambda}}}(G(x,\bar{y})-G(\bar{x},\bar{y}))\tilde{f}(u(\bar{y}))dy
\end{align*}
By previous argument, we know that either $f(u(y))>\tilde{f}(u(\bar{y}))$ or $\tilde{f}(u(\bar{y}))>0$ in $\mathcal{O}_{\bar{\lambda}}\subset\Sigma_\lambda$. In either case, we have
$$u(x)-u(\bar{x})>\int_{\Sigma_{\bar{\lambda}}}(G(x,y)-G(\bar{x},y)+G(x,\bar{y})-G(\bar{x},\bar{y}))\tilde{f}(u(\bar{y}))dy\geq 0,$$ due to Lemma \ref{max G}.
Thus, $u(x)>u(\bar{x})$ for all $x\in\Sigma_{\bar{\lambda}}$. Then by compactness, there exists $0<\gamma_1<\gamma$ such that $u(x)>u(\bar{x})$
for all $l\in(\bar{\lambda}-\gamma_1,\bar{\lambda}]$. which contradicts the choices of $\bar{\lambda}$, together with (\ref{cty}).
Therefore, since $\Lambda=(0, R)$, we have $u(\bar{x}_1,x_2,\cdots,x_n)\geq u(x_1,x_2,\cdots,x_n)$ if $x_1\geq 0$. We conclude that $u$ is radially symmetric by similar argument from the proof of the last theorem.

\subsection{Proof of Theorem \ref{thm3}}
If $u$ is a positive solution to (\ref{P_k}), $u$ is radial symmetric with respect to the origin.
Let $$v=\left(\frac{1-|x|^2}{2}\right)^{k-\frac{n}{2}}u,$$ then $v$ is a positive solution to (\ref{eq1.7}) and is again symmetric with respect to the origin on $B(0,R^\prime)$.
We are left to show that the boundary condition $\nabla^\alpha v|_{\partial B^\prime}=0$ holds. We prove this by induction. For $\alpha=1$, it is clear that on the boundary,
\begin{align*}
    \nabla_{\hs}u & = \left(\frac{1-|x|^2}{2}\nabla\right)\left(\frac{1-|x|^2}{2}\right) ^{\frac{n}{2}-k}v\\
    &=\left(\frac{1-|x|^2}{2}\right)\left(\frac{1-|x|^2}{2}\nabla v+\nabla\left(\frac{1-|x|^2}{2}\right)^{\frac{n}{2}-k}v\right)\\
    &=\left(\frac{1-|x|^2}{2}\right)^2\nabla v=0 \text{ on } \partial B^\prime.
\end{align*}
Thus, $\nabla v|_{\partial B^\prime}=0$. Suppose this is also for true for $1\leq\alpha<k-2$, then we have that on $\partial B^\prime$,
\begin{align*}
    \nabla^{\alpha+1}_{\hs}u & = \left(\frac{1-|x|^2}{2}\nabla\right)\left(\frac{1-|x|^2}{2}\nabla\right)^{\alpha}\left[\left(\frac{1-|x|^2}{2}\right) ^{\frac{n}{2}-k}v\right]\\
    &=\left(\frac{1-|x|^2}{2}\nabla\right) \left(\frac{1-|x|^2}{2}\right)^{\alpha}\nabla^{\alpha}v+\sum_{i=0}^{\alpha}\nabla_{\hs}^i u\\
    &=\left(\frac{1-|x|^2}{2}\right)^{\alpha+1}\nabla^{\alpha+1}v+\nabla^\alpha v\left(\frac{1-|x|^2}{2}\nabla\right)\left(\frac{1-|x|^2}{2}\right)^{\alpha}\\
    &=\left(\frac{1-|x|^2}{2}\right)^{\alpha+1}\nabla^{\alpha+1}v=0.
\end{align*}
Therefore, $\nabla^{\alpha}v=0$ for all $|\alpha|\leq k-1$.

\subsection{Proof of Corollary \ref{cor1.4}}
It is known that $P_k=|x_n|^{k+\frac{n}{2}}(-\Delta)^k\left(|x_n|^{k-\frac{n}{2}}u\right)$ on $\mathbb{R}^n_+$. Thus, let $w=|x_n|^{k-\frac{n}{2}}u$,
$P_ku=f(u)$ gives
$$(-\Delta)^kw = \frac{1}{|x_n|^{k+\frac{n}{2}}}f(|x_n|^{\frac{n}{2}-k}w)\; \text{ on } \halfrn.$$
From Theorem \ref{thm1}, we know that $u$ is radial symmetric.
That is, $u(x)=u(\bar{x})$, if $\rho(x,p)=\rho(\bar{x},p)$, the geodesic distance to $p$.
\\~\\
Let $p=(p^\prime,p_n)\in\mathbb{R}^n_+$, where $p^\prime=(p_1,\cdots,p_{n-1})$ with $p_1=\lambda$. Without loss of generality, we choose $x_1=\lambda$ to be the direction which is parallel to $x_n$-axis
and let $T$ be the reflection map with respect to $x_1=\lambda$.
For any $x=(x^\prime,x_n)$, let  $\bar{x}=T(x)=(\bar{x}^\prime,\bar{x}_n)$.
We first claim $\rho(x,p)=\rho(\bar{x},p)$, where $\rho(\cdot,p)$ is the distance to $p$. The claim is true by observing that $T$ is an isometry so that
 $d(x,p)=\rho(T(x),T(p))=\rho(\bar{x},p)$. Thus, $w(x)=w(\bar{x})$. The solutions are axially symmetric about some line parallel to $x_n$-axis.

\subsection{Proof of Theorem \ref{thm1.5} and \ref{thm1.6}}
We follow the similar notations in the proof of Theorem \ref{thm1} and let
$U=\bn\cap\{x_1=0\}$ and  $U_\lambda=A_\lambda(U), \Sigma_\lambda=\cup_{s<\lambda}U_s$.
For any $x\in\Sigma_\lambda$, denote $\bar{x}=I_\lambda(x)$ and $u_\lambda(x)=u(\bar{x})$. The proof will again be complete once we show that $\lambda=0$.
Now assume by contradiction that $\lambda>0$, we get

\begin{align*}
    &u(x)-u_\lambda(x)
    =\int_{\Sigma_\lambda} \left(G_\alpha(x,y)-G_\alpha(\bar{x},y)\right)\left(f(u)-f(u_\lambda)\right)dV_y
    \end{align*}
    In $\Sigma_\lambda$, we denote $\Sigma_\lambda^-=\{x\in\Sigma_\lambda: u_\lambda(x)>u(x)\}$, and we would like to show that $\Sigma_\lambda^-$ is of measure zero.
    In fact, by (\ref{k_est}) and the fact that $\rho\sim\sinh\rho$ for small $\rho$, we have
    \begin{align*}
    u_\lambda(x)-u(x)&\leq \int_{\Sigma_\lambda^-}G_\alpha(x,y)\left(f(u_\lambda)-f(u)\right)dV_y\\
    &\leq C\int_{\Sigma_\lambda^-}\left(\frac{1}{2\sinh \frac{\rho}{2}}\right)^{n-\alpha}\left(f(u_\lambda)-f(u)\right)dV_y.
    \end{align*}
    Thus, the Hardy-Littlewood-Sobolev inequality on $\hs$ can again be appied to here, together with mean value theorem and H\"older's inequality, we have
    \begin{align*}
    \|u-u_\lambda\|_{L^q(\Sigma_\lambda)}&\leq C\|f(u)-f(u_\lambda)\|_{\frac{nq}{n+\alpha q}}\\
    &=C\left(\int_{\Sigma_\lambda^-}\left(f^\prime(\xi)(u(y)-u_\lambda(y))\right)^{\frac{nq}{n+\alpha q}}\right)^{\frac{1}{q}+\frac{\alpha}{n}}\\
    %&\leq C\left(\int_{\Sigma_t^-}f^\prime(\xi)^\frac{n}{2k}\right)^{\frac{2k}{n}}\left(\int_{\Sigma_t^-}\left(u(y)-u_t(y)\right)^q\right)^\frac{1}{q}\\%
    &\leq C\|f^\prime(\xi)\|_{L^{\frac{n}{\alpha}}(\Sigma_\lambda^-)}\|u-u_\lambda\|_{L^q(\Sigma_\lambda^-)},
    \end{align*}
    for some $\frac{n}{n-\alpha}<q<\infty$.
    Here $\xi=\theta u + (1-\theta)u_\lambda$ where $0<\theta<1$.
    Since either (i) or (ii) holds for $f(u)$, we are able to pick $\lambda$ large enough so that $C\|f^\prime(\xi)\|_{\frac{n}{\alpha}}<1$.
    This implies that for large enough $\lambda, \|u-u_\lambda\|_{L^q(\Sigma_\lambda^-)}=0$, which shows that $\Sigma_\lambda^-$ is of measure zero.
    \\~\\
    Now we shift $U_\lambda$ as long as $u\geq u_\lambda$ in $\Sigma_\lambda$. Suppose that there exists such $\bar{\lambda}$ that $u(x) > u_{\bar{\lambda}}(x)$ on $\Sigma_{\bar{\lambda}}$.
    We deduce again by a standard compactness argument that there exist $\bar{\lambda}-\varepsilon<\lambda\leq \bar{\lambda}$ such that
    $$\|u-u_\lambda\|_{L^q(\Sigma_\lambda)}\leq C\|f^\prime(\xi)\|_{L^{\frac{n}{\alpha}}(\Sigma_\lambda^-)}\|u-u_\lambda\|_{L^q(\Sigma_\lambda^-)}.$$
    When $\varepsilon$ is small, $\Sigma_\lambda^-$ is close to zero so that $\|u_t-u\|_{L^q(\Sigma_\lambda^-)}=0$.
    This imples we can keep moving the plane $U_{\lambda_0}$ until we reach the origin. Now we see $u(x)\leq u(\bar{x})$ with respect to $\Sigma_0$.
    A similar argument shows $u(\bar{x})\geq u(x)$ by rotation. Consequently, $u$ is symmetric with respect to every hyperplane containing the origin, namely, is radial symmetric.
\\~\\
The proof of Theorem \ref{thm1.6} is similar, the essential part is also to show that $\Sigma_\lambda^-$ has measure zero so that the plane can be moved from the infinity. To achieve this, we notice that
\begin{align*}
    u_\lambda(x)-u(x)&\leq \int_{\Sigma_\lambda^-}\mathcal{G}(x,y)\left(f(u_\lambda)-f(u)\right)dV_y\\
    &\leq C\int_{\Sigma_\lambda^-}\left(\frac{1}{2\sinh \frac{\rho}{2}}\right)^{n-S_l}\left(f(u_\lambda)-f(u)\right)dV_y.
    \end{align*}
Therefore, by the Hardy-Littlewood-Sobolev inequality on $\hs$, mean value theorem and H\"older's inequality, we have
\begin{align*}
\|u-u_\lambda\|_{L^q(\Sigma_\lambda)}&\leq C\|f(u)-f(u_\lambda)\|_{\frac{nq}{n+S_lq}}\\
&=C\left(\int_{\Sigma_\lambda^-}\left(f^\prime(\xi)(u(y)-u_\lambda(y))\right)^{\frac{nq}{n+S_lq}}\right)^{\frac{1}{q}+\frac{S_l}{n}}\\
%&\leq C\left(\int_{\Sigma_t^-}f^\prime(\xi)^\frac{n}{2k}\right)^{\frac{2k}{n}}\left(\int_{\Sigma_t^-}\left(u(y)-u_t(y)\right)^q\right)^\frac{1}{q}\\%
&\leq C\|f^\prime(\xi)\|_{L^{\frac{n}{S_l}}(\Sigma_\lambda^-)}\|u-u_\lambda\|_{L^q(\Sigma_\lambda^-)},
\end{align*}
for some $\frac{n}{n-S_l}<q<\infty$. The same argument shows that $\Sigma_\lambda^-$ has measure zero.


\begin{thebibliography}{99}

    \bibitem{Ahlfors} Ahlfors, L. V. \textit{M\"obius trasnformations in several dimensions, Ordway Professorship Lectures in Mathematics}, University of Minnesota, School of Mathematics, Minneapolis, MN. 1981.

    \bibitem{ADG} Almeida, L; Damascelli, L; Ge, Y. \textit{A few symmetry results for nonlinear elliptic PDE on noncompact manifolds}, Ann. Inst. H. Poinc\'are Anal. Non Lin\'eaire, 19 (2002) 313-342.

    \bibitem{A-J}  Anker, J.P.;  Ji, L. \textit{Heat kernel and Green function estimates on noncompact symmetric spaces}, Geom. Funct. Anal. 9 (1999) 1035-1091.

    \bibitem{Yamabe1} Aubin, T. \textit{Équations différentielles non linéaires et problème de Yamabe concernant la courbure scalaire}, J. Math. Pures Appl., 55 (1976) 269-296.

    \bibitem{Beckner} Beckner, W. \textit{On Lie groups and hyperbolic symmetry--from Kunze-Stein phenomena to Riesz potentials}, Nonlinear Anal. 126 (2015) 394-414.

    \bibitem{Benguria} Benguria, S. \textit{The solution gap of the Br\'ezis-Nirenberg problem on the hyperbolic space}, Monatsh. Math., 181 (2016) 537-559.

    \bibitem{BGW} Berchio, E; Gazzola, F; Weth, T. \textit{Radial symmetry of positive solutions to nonlinear polyharmonic Dirichlet problems}. J. Reine Angew. Math., 2008 (2008) 165-183.

    \bibitem{Bere} Berestycki, H; Nirenberg, L. \textit{On the method of moving planes and the sliding method}, Bol. Soc. Bras. Mat., 22 (1991) 1-37.

    \bibitem{BernisGrunau1} Bernis, F; Grunau, H-C. \textit{Critical exponents and multiple critical dimensions for polyharmonic operators}, J. Differential Equations, 117 (1995) 469-486.

    \bibitem{Bo} Boggio, T. \textit{Sulle funzioni di Green d'ordine m}, Rend. Circ. Mat. Palermo, 20 (1905) 07-135.

    \bibitem{B-N} Br\'ezis, H; Nirenberg, L. \textit{Positive solutions of nonlinear elliptic equations involving critical sobolev exponents}, Commun. Pure Appl. Anal., 36 (1983) 437-477.

    \bibitem{Caff} Caffarelli, L; Gidas, B; Spruck, J. \textit{Asymptotic symmetry and local behavior of semilinear elliptic equations with critical Sobolev growth}, Comm. Pure Appl. Math., 42(3) (1989) 271-297.

    \bibitem{ChenLi1} Chen, W.; Li, C. \textit{Methods on nonlinear elliptic equations}, AIMS Series on Differential Equations \& Dynamical Systems, American Institute of Mathematical Sciences (AIMS), Springfield, MO, 2010. xii+299 pp.

    \bibitem{CLO} Chen, W; Li, C; Ou, B. \textit{Classification of solutions for an integral equation}, Comm. Pure Appl. Math., 59(3) (2006) 330-343.

    
    \bibitem{D-M} Davies, E. B.; Mandouvalos, N. \textit{Heat kernel bounds on hyperbolic space and Kleinian groups}, Proc. London Math. Soc. (3) 57 (1998) 182-208.

    \bibitem{EdmundsFortunatoJannelli1} Edmunds, D; Fortunato, D; Jannelli,E. \textit{Critical exponents, critical dimensions and the biharmonic operator}, Arch. Rational Mech. Anal., 112(3) (1990) 269-289.

    \bibitem{EM} Erdélyi, A; Magnus, W; Oberhettinger, F; Tricomi,  F.G. \textit{Higher Transcendental Functions}, vol. I, Based, in part, on notes left by Harry Bateman, McGraw-Hill Book Company, Inc., New York-Toronto-London, (1953).

    \bibitem{FeffermanGr} Fefferman, C; Graham, C. R.  The ambient metric, Ann. of Math. Stud., Number 178, Princeton University Press, (2012).

    \bibitem{FeffermanGr1}  Fefferman, C; Graham, C. R. Juhl's formulae for GJMS operators and Q-curvatures. J. Amer. Math. Soc. 26 (2013), no. 4, 1191-1207.

    \bibitem{Gazzola} Gazzola, F. Critical growth problems for polyharmonic operators. Proc. Roy. Soc. Edinburgh Sect. A 128 (1998), no. 2, 251-263.

    \bibitem{GNN} Gidas, B; Ni, W; Nirenberg, L. \textit{Symmetry and related properties via the maximum principle}, Comm. Math. Phys., 68(3) (1979) 209-243.

    \bibitem{GJMS1} Gover, A. R.  \textit{Laplacian operators and Q-curvature on conformally Einstein manifolds}, Math. Ann. 336 (2006) 311-334.

    \bibitem{GJMS2} Graham, C. R; Jenne, R; Mason, L. J; Sparling, G. A. J.\textit{Conformally Invariant Powers of the Laplacian, I: Existence}, J. London Math. Soc. (2) 46 (1992) no. 3, 557-565.

    \bibitem{G-N} Grigor'yan A; Noguchi, M. \text{The hear kernelon hyperbolic space}, Bull. London Math. Soc. 30 (1998), 30, 643-650.

    \bibitem{Grunau} Grunau, H.-C. \textit{Positive solutions to semilinear polyharmonic Dirichlet problems involving critical Sobolev exponents}, Calc. Var. Partial Differential Equations, 3.2 (1995) 243-252.

    \bibitem{H-F1} Helgason, S. \textit{Groups and geometric analysis}, Mathematical Surveys and Monographs, vol. 83, 1984.

    \bibitem{H-F2} Helgason, S. \textit{Geometric analysis on symmetric spaces}, Mathematical Surveys and Monographs, vol: 39, 2008.

    \bibitem{GJMS3} Juhl, A. \textit{Explicit Formulas for GJMS-Operators and Q-Curvatures}, Geom. Funct. Anal. 23 (2013) 1278-1370.

    \bibitem{K-P1} Kumaresan, S; Prajapat, J. V. \textit{Analogue of Gidas-Ni-Nirenberg result in hyperbolic space and sphere}, Rend. Istit. Mat. Univ. Trieste 30 (1998), no. 1-2, 107-112.

    \bibitem{K-P2} Kumaresan, S; Prajapat, J. V. \textit{Serrin's result for hyperbolic space and sphere}, Duke. Math. J., Vol. 91 (1998) 17-28.

    \bibitem{LeeParker1} Lee, J. M; Parker; T. H. \textit{The Yamabe problem}, Bull. Amer. Math. Soc. (N.S.), 17(1) (1987) 37-91.

    \bibitem{Lhq} Li, H.Q. \textit{Fonction maximale centr\'ee de Hardy-Littlewood sur les espaces hyperboliques}, Ark. Mat. 50 (2012), no. 2, 359-378.

    \bibitem{LLY2} Li, J; Lu, G; Yang, Q. \textit{Sharp Adam and Hardy-Adams inequalities of any fractional order on hyperbolic spaces of all dimensions}, Trans. Amer. Math. Soc. 373 (2020), no. 5, 3483-3513.

    \bibitem{LLY}  Li, J; Lu, G; Yang, Q. \textit{Higher order Brezis-Nirenberg problem on hyperbolic spaces: Existence, nonexistence and symmetry of solutions}, Adv. Math. 399 (2022), Paper No. 108259, 39 pp.

    \bibitem{Li-Yau} Li, P; Yau, S.T. \textit{On the parabolic kernel of the Schrödinger operator}, Acta Math. 156 (1986) 153-201.

    \bibitem{Liu} Liu, G. \textit{Sharp higher-order Sobolev inequalities in the hyperbolic space $\mathbb{H}^n$}, Calc. Var. Partial Differential Equations, 47(3-4) (2013) 567-588.

    \bibitem{L-Y} Lu, G; Yang, Q. \textit{Paneitz operators on hyperbolic spaces and high order Hardy-Sobolev-Maz'ya inequalities on half spaces}, Amer. J. Math., vol. 141 no. 6, (2019), p. 1777-1816.

    \bibitem{LY} Lu, G; Yang, Q. \textit{Green's functions of Paneitz and GJMS operators on hyperbolic spaces and sharp Hardy-Sobolev-Maz'ya inequaities on half spaces}, Adv. Math. 398 (2022), Paper No. 108156, 42 pp.

    \bibitem{L-Z} Lu, G; Zhu, J. \textit{Axial symmetry and regularity of solutions to an integral equation in a half-space}, Pacific J. Math. 253 (2011), no. 2, 455-473.

    \bibitem{Mancini-Sandeep} Mancini, G; Sandeep, K. \textit{On a semilinear elliptic equation in $\mathbb{H}^n$}, Ann. Sc. Norm. Super. Pisa Cl. Sci. (5), 7(4) (2008) 635-671.

    \bibitem{HM} Matsumoto, H. \textit{Closed form formulae for the heat kernels and the Green functions for the Laplacians on the symmetric spaces of rank one}, Bull. Sci. Math. 125 (2001) 553-581.

    \bibitem{PucciSerrin1} Pucci, P; Serrin, J. \textit{Critical exponents and critical dimensions for polyharmonic operators}, J. Math. Pures Appl. (9), 69(1) (1990) 55--83.

    \bibitem{Ratcliffe} Ratcliffe, J. G. \textit{Foundations of Hyperbolic Manifolds}, Springer-Verlag New York 2006.

    \bibitem{Yamabe2} Schoen, R. \textit{Conformal deformation of a Riemannian metric to constant scalar curvature},  J. Differential Geom. 20 (1984), no. 2, 479-495.

    \bibitem{serrin} Serrin, J. \textit{A symmetry problem in potential theory}, Arch. Rat. Mech. 43 (1971) 304-318.

    \bibitem{Stapelkamp} Stapelkamp, S. \textit{The {B}r\'{e}zis-{N}irenberg problem on $\hs$. Existence and uniqueness of solutions}, J. Elliptic Parabol. Equ., (2002) 283-290.

    \bibitem{Strichartz} Strichartz, R. S. \textit{Analysis of the Laplacian on the complete Riemannian manifold}, J. Funct. Anal., 52, 1 (1983) 48-79.

    \bibitem{Talenti1} Talenti, G. \textit{Best constant in {S}obolev inequality}, Ann. Mat. Pura Appl. (4), 110 (1976) 353--372.

    \bibitem{Yamabe3} Trudinger, N.S. \textit{Remarks concerning the conformal deformation of Riemannian structures on compact manifolds}, Ann. Scuola Norm. Sup. Pisa (3), 22 (1968) 265-274.

    \bibitem{W-X} Wei, J; Xu, X. \textit{Classification of solutions of higher order conformally invariant equations}, Math. Ann. 313 (1999) 207-228.

    \bibitem{Yamabe} Yamabe, H. \textit{On a deformation of Riemannian structures on compact manifolds}, Osaka Math. J. 12 (1960) 21-37.
    \end{thebibliography}
\end{document}